\pdfoutput=1
\documentclass[a4paper,10pt]{amsart}

\usepackage[utf8]{inputenc}
\usepackage{url}

\usepackage[expansion=false]{microtype}
\usepackage{amssymb,amsmath,amsopn,amsxtra,amsthm,amsfonts}
\usepackage{xr}
\usepackage[mathcal]{euscript}
\usepackage{mathalfa}

\usepackage[english, status=final, nomargin, inline]{fixme}
\fxusetheme{color}

\usepackage[pdfusetitle,unicode,hidelinks]{hyperref}

\makeatletter
\let\c@equation\c@subsubsection

\makeatother

\newtheorem{cor}[subsubsection]{Corollary}
\newtheorem{lem}[subsubsection]{Lemma}
\newtheorem{prop}[subsubsection]{Proposition}

\newtheorem{thm}[subsubsection]{Theorem}

\newtheorem{defn}[subsubsection]{Definition}

\newtheorem*{claim*}{Claim}

\theoremstyle{definition}

\theoremstyle{remark}
\newtheorem{rem}[subsubsection]{Remark}

\newtheorem{ex}[subsubsection]{Example}

\newcommand{\secref}[1]{Section~\ref{#1}}
\newcommand{\ssecref}[1]{\sectsign\ref{#1}}
\newcommand{\sssecref}[1]{\ref{#1}}

\renewcommand{\eqref}[1]{(\ref{#1})}

\usepackage{tikz}
\usetikzlibrary{matrix}
\usepackage{tikz-cd}

\usepackage{xspace}

\usepackage{xifthen}

\usepackage{xparse}

\newcommand{\nc}{\newcommand}
\nc{\renc}{\renewcommand}
\nc{\ssec}{\subsection}
\nc{\sssec}{\subsubsection}
\nc{\on}{\operatorname}
\nc{\term}[1]{#1\xspace}

\DeclareMathSymbol{A}{\mathalpha}{operators}{`A}
\DeclareMathSymbol{B}{\mathalpha}{operators}{`B}
\DeclareMathSymbol{C}{\mathalpha}{operators}{`C}
\DeclareMathSymbol{D}{\mathalpha}{operators}{`D}
\DeclareMathSymbol{E}{\mathalpha}{operators}{`E}
\DeclareMathSymbol{F}{\mathalpha}{operators}{`F}
\DeclareMathSymbol{G}{\mathalpha}{operators}{`G}
\DeclareMathSymbol{H}{\mathalpha}{operators}{`H}
\DeclareMathSymbol{I}{\mathalpha}{operators}{`I}
\DeclareMathSymbol{J}{\mathalpha}{operators}{`J}
\DeclareMathSymbol{K}{\mathalpha}{operators}{`K}
\DeclareMathSymbol{L}{\mathalpha}{operators}{`L}
\DeclareMathSymbol{M}{\mathalpha}{operators}{`M}
\DeclareMathSymbol{N}{\mathalpha}{operators}{`N}
\DeclareMathSymbol{O}{\mathalpha}{operators}{`O}
\DeclareMathSymbol{P}{\mathalpha}{operators}{`P}
\DeclareMathSymbol{Q}{\mathalpha}{operators}{`Q}
\DeclareMathSymbol{R}{\mathalpha}{operators}{`R}
\DeclareMathSymbol{S}{\mathalpha}{operators}{`S}
\DeclareMathSymbol{T}{\mathalpha}{operators}{`T}
\DeclareMathSymbol{U}{\mathalpha}{operators}{`U}
\DeclareMathSymbol{V}{\mathalpha}{operators}{`V}
\DeclareMathSymbol{W}{\mathalpha}{operators}{`W}
\DeclareMathSymbol{X}{\mathalpha}{operators}{`X}
\DeclareMathSymbol{Y}{\mathalpha}{operators}{`Y}
\DeclareMathSymbol{Z}{\mathalpha}{operators}{`Z}

\nc{\sA}{\ensuremath{\mathcal{A}}\xspace}
\nc{\sB}{\ensuremath{\mathcal{B}}\xspace}
\nc{\sC}{\ensuremath{\mathcal{C}}\xspace}
\nc{\sD}{\ensuremath{\mathcal{D}}\xspace}
\nc{\sE}{\ensuremath{\mathcal{E}}\xspace}
\nc{\sF}{\ensuremath{\mathcal{F}}\xspace}
\nc{\sG}{\ensuremath{\mathcal{G}}\xspace}
\nc{\sH}{\ensuremath{\mathcal{H}}\xspace}
\nc{\sI}{\ensuremath{\mathcal{I}}\xspace}
\nc{\sJ}{\ensuremath{\mathcal{J}}\xspace}
\nc{\sK}{\ensuremath{\mathcal{K}}\xspace}
\nc{\sL}{\ensuremath{\mathcal{L}}\xspace}
\nc{\sM}{\ensuremath{\mathcal{M}}\xspace}
\nc{\sN}{\ensuremath{\mathcal{N}}\xspace}
\nc{\sO}{\ensuremath{\mathcal{O}}\xspace}
\nc{\sP}{\ensuremath{\mathcal{P}}\xspace}
\nc{\sQ}{\ensuremath{\mathcal{Q}}\xspace}
\nc{\sR}{\ensuremath{\mathcal{R}}\xspace}
\nc{\sS}{\ensuremath{\mathcal{S}}\xspace}
\nc{\sT}{\ensuremath{\mathcal{T}}\xspace}
\nc{\sU}{\ensuremath{\mathcal{U}}\xspace}
\nc{\sV}{\ensuremath{\mathcal{V}}\xspace}
\nc{\sW}{\ensuremath{\mathcal{W}}\xspace}
\nc{\sX}{\ensuremath{\mathcal{X}}\xspace}
\nc{\sY}{\ensuremath{\mathcal{Y}}\xspace}
\nc{\sZ}{\ensuremath{\mathcal{Z}}\xspace}

\nc{\bA}{\ensuremath{\mathbf{A}}\xspace}
\nc{\bB}{\ensuremath{\mathbf{B}}\xspace}
\nc{\bC}{\ensuremath{\mathbf{C}}\xspace}
\nc{\bD}{\ensuremath{\mathbf{D}}\xspace}
\nc{\bE}{\ensuremath{\mathbf{E}}\xspace}
\nc{\bF}{\ensuremath{\mathbf{F}}\xspace}
\nc{\bG}{\ensuremath{\mathbf{G}}\xspace}
\nc{\bH}{\ensuremath{\mathbf{H}}\xspace}
\nc{\bI}{\ensuremath{\mathbf{I}}\xspace}
\nc{\bJ}{\ensuremath{\mathbf{J}}\xspace}
\nc{\bK}{\ensuremath{\mathbf{K}}\xspace}
\nc{\bL}{\ensuremath{\mathbf{L}}\xspace}
\nc{\bM}{\ensuremath{\mathbf{M}}\xspace}
\nc{\bN}{\ensuremath{\mathbf{N}}\xspace}
\nc{\bO}{\ensuremath{\mathbf{O}}\xspace}
\nc{\bP}{\ensuremath{\mathbf{P}}\xspace}
\nc{\bQ}{\ensuremath{\mathbf{Q}}\xspace}
\nc{\bR}{\ensuremath{\mathbf{R}}\xspace}
\nc{\bS}{\ensuremath{\mathbf{S}}\xspace}
\nc{\bT}{\ensuremath{\mathbf{T}}\xspace}
\nc{\bU}{\ensuremath{\mathbf{U}}\xspace}
\nc{\bV}{\ensuremath{\mathbf{V}}\xspace}
\nc{\bW}{\ensuremath{\mathbf{W}}\xspace}
\nc{\bX}{\ensuremath{\mathbf{X}}\xspace}
\nc{\bY}{\ensuremath{\mathbf{Y}}\xspace}
\nc{\bZ}{\ensuremath{\mathbf{Z}}\xspace}

\nc{\dA}{\ensuremath{\mathds{A}}\xspace}
\nc{\dB}{\ensuremath{\mathds{B}}\xspace}
\nc{\dC}{\ensuremath{\mathds{C}}\xspace}
\nc{\dD}{\ensuremath{\mathds{D}}\xspace}
\nc{\dE}{\ensuremath{\mathds{E}}\xspace}
\nc{\dF}{\ensuremath{\mathds{F}}\xspace}
\nc{\dG}{\ensuremath{\mathds{G}}\xspace}
\nc{\dH}{\ensuremath{\mathds{H}}\xspace}
\nc{\dI}{\ensuremath{\mathds{I}}\xspace}
\nc{\dJ}{\ensuremath{\mathds{J}}\xspace}
\nc{\dK}{\ensuremath{\mathds{K}}\xspace}
\nc{\dL}{\ensuremath{\mathds{L}}\xspace}
\nc{\dM}{\ensuremath{\mathds{M}}\xspace}
\nc{\dN}{\ensuremath{\mathds{N}}\xspace}
\nc{\dO}{\ensuremath{\mathds{O}}\xspace}
\nc{\dP}{\ensuremath{\mathds{P}}\xspace}
\nc{\dQ}{\ensuremath{\mathds{Q}}\xspace}
\nc{\dR}{\ensuremath{\mathds{R}}\xspace}
\nc{\dS}{\ensuremath{\mathds{S}}\xspace}
\nc{\dT}{\ensuremath{\mathds{T}}\xspace}
\nc{\dU}{\ensuremath{\mathds{U}}\xspace}
\nc{\dV}{\ensuremath{\mathds{V}}\xspace}
\nc{\dW}{\ensuremath{\mathds{W}}\xspace}
\nc{\dX}{\ensuremath{\mathds{X}}\xspace}
\nc{\dY}{\ensuremath{\mathds{Y}}\xspace}
\nc{\dZ}{\ensuremath{\mathds{Z}}\xspace}

\nc{\bbA}{\ensuremath{\mathbb{A}}\xspace}
\nc{\bbB}{\ensuremath{\mathbb{B}}\xspace}
\nc{\bbC}{\ensuremath{\mathbb{C}}\xspace}
\nc{\bbD}{\ensuremath{\mathbb{D}}\xspace}
\nc{\bbE}{\ensuremath{\mathbb{E}}\xspace}
\nc{\bbF}{\ensuremath{\mathbb{F}}\xspace}
\nc{\bbG}{\ensuremath{\mathbb{G}}\xspace}
\nc{\bbH}{\ensuremath{\mathbb{H}}\xspace}
\nc{\bbI}{\ensuremath{\mathbb{I}}\xspace}
\nc{\bbJ}{\ensuremath{\mathbb{J}}\xspace}
\nc{\bbK}{\ensuremath{\mathbb{K}}\xspace}
\nc{\bbL}{\ensuremath{\mathbb{L}}\xspace}
\nc{\bbM}{\ensuremath{\mathbb{M}}\xspace}
\nc{\bbN}{\ensuremath{\mathbb{N}}\xspace}
\nc{\bbO}{\ensuremath{\mathbb{O}}\xspace}
\nc{\bbP}{\ensuremath{\mathbb{P}}\xspace}
\nc{\bbQ}{\ensuremath{\mathbb{Q}}\xspace}
\nc{\bbR}{\ensuremath{\mathbb{R}}\xspace}
\nc{\bbS}{\ensuremath{\mathbb{S}}\xspace}
\nc{\bbT}{\ensuremath{\mathbb{T}}\xspace}
\nc{\bbU}{\ensuremath{\mathbb{U}}\xspace}
\nc{\bbV}{\ensuremath{\mathbb{V}}\xspace}
\nc{\bbW}{\ensuremath{\mathbb{W}}\xspace}
\nc{\bbX}{\ensuremath{\mathbb{X}}\xspace}
\nc{\bbY}{\ensuremath{\mathbb{Y}}\xspace}
\nc{\bbZ}{\ensuremath{\mathbb{Z}}\xspace}

\nc{\mrm}[1]{\ensuremath{\mathrm{#1}}\xspace}
\nc{\mbf}[1]{\ensuremath{\mathbf{#1}}\xspace}
\nc{\mcal}[1]{\ensuremath{\mathcal{#1}}\xspace}
\nc{\msc}[1]{\ensuremath{\mathscr{#1}}\xspace}

\renc{\bar}[1]{\overline{#1}}

\let\sectsign\S
\let\S\relax

\nc{\sub}{\subset}
\nc{\too}{\longrightarrow}
\nc{\hook}{\hookrightarrow}
\nc*{\hooklongrightarrow}{\ensuremath{\lhook\joinrel\relbar\joinrel\rightarrow}}
\nc{\hooklong}{\hooklongrightarrow}
\nc{\twoheadlongrightarrow}{\relbar\joinrel\twoheadrightarrow}
\nc{\shiso}{\approx}
\nc{\isoto}{\xrightarrow{\sim}}
\nc{\isofrom}{\xleftarrow{\sim}}
\renc{\ge}{\geqslant}
\renc{\le}{\leqslant}
\renc{\geq}{\geqslant}
\renc{\leq}{\leqslant}

\nc{\id}{\mathrm{id}}

\DeclareMathOperator{\rk}{\mathrm{rk}}
\DeclareMathOperator{\Hom}{\mathrm{Hom}}
\nc{\uHom}{\underline{\smash{\Hom}}}
\DeclareMathOperator{\Maps}{\mathrm{Maps}}

\DeclareMathOperator{\End}{\mathrm{End}}
\DeclareMathOperator{\Sym}{\mathrm{Sym}}
\nc{\Pre}{\mathrm{PSh}{}}
\nc{\uEnd}{\underline{\smash{\End}}}
\DeclareMathOperator{\codim}{\mathrm{codim}}
\renc{\lim}{\operatorname*{lim}}
\nc{\colim}{\operatorname*{colim}}
\nc{\Cofib}{\on{Cofib}}
\nc{\Fib}{\on{Fib}}
\nc{\initial}{\varnothing}
\nc{\op}{\mathrm{op}}

\DeclareMathOperator*{\fibprod}{\times}

\renc{\coprod}{\sqcup}

\nc{\bDelta}{\mbf{\Delta}}
\nc{\DM}{\mbf{DM}}
\nc{\eff}{\mathrm{eff}}
\nc{\veff}{\mathrm{veff}}
\nc{\cyc}{{\mrm{cyc}}}
\nc{\corr}{{\on{corr}}}
\nc{\ft}{\mrm{ft}}
\nc{\flf}{\mrm{flf}}
\nc{\fet}{{\mrm{f\acute et}}}
\nc{\fsyn}{{\mrm{fsyn}}}
\nc{\syn}{{\mrm{syn}}}
\nc{\lci}{{\mrm{lci}}}
\nc{\Perf}{\mbf{Perf}}
\nc{\perf}{\on{perf}}
\nc{\oblv}{\on{oblv}}
\nc{\exact}{\on{exact}}

\nc{\F}{{\on{F}}}
\nc{\clopen}{{\mrm{clopen}}}
\nc{\B}{\mrm{B}}
\nc{\D}{\mrm{D}}
\nc{\Fin}{\on{Fin}}
\nc{\Cut}{\on{Cut}}
\nc{\Cart}{\on{Cart}}
\nc{\pairs}{\mathsf{pairs}}
\nc{\Pairs}{\mathrm{Pair}}
\nc{\Trip}{\mathrm{Trip}}
\nc{\Lab}{\mathrm{Lab}}
\nc{\SL}{\mathrm{SL}}
\nc{\coCart}{\mathrm{coCart}}
\nc{\RKE}{\mathrm{RKE}}
\nc{\strict}{\mathrm{strict}}
\nc{\Emb}{\mathrm{Emb}}
\nc{\Split}{\mathrm{Split}}
\nc{\Set}{\mathrm{Set}}
\nc{\sSets}{\mathrm{sSets}}
\nc{\pb}{\mathrm{pb}}
\nc{\fib}{\mathrm{fib}}
\nc{\diff}{\mrm{diff}}
\nc{\gp}{\mrm{gp}}
\DeclareMathOperator{\chr}{char}
\nc{\mgp}{\mrm{mot-gp}}
\nc{\FSyn}{\mrm{FSyn}}
\nc{\FEt}{\mrm{FEt}}
\nc{\Spc}{\mrm{Spc}}
\nc{\Ob}{\mrm{Ob}}
\nc{\Spt}{\mrm{Spt}}
\nc{\T}{\bT}
\nc{\suspinf}{\Sigma^\infty}
\nc{\h}{\mrm{h}}
\nc{\uhom}{\underline{\mathrm{Hom}}}
\nc{\umap}{\underline{\mathrm{Maps}}}
\renc{\H}{\bH}
\nc{\Einfty}{{\sE_\infty}}
\nc{\Eone}{{\sE_1}}
\nc{\Stab}{\mrm{Stab}}
\nc{\lax}{{\mrm{lax}}}
\nc{\cocart}{{\mrm{cocart}}}
\nc{\Sch}{\on{Sch}}
\nc{\Fr}{\on{Fr}}
\nc{\A}{\mathbf{A}}
\nc{\N}{\mathbf{N}}
\nc{\Z}{\mathbf{Z}}
\nc{\Q}{\mathbf{Q}}
\nc{\Oo}{\mathcal{O}} 
\nc{\Ll}{\mathcal{L}} 
\nc{\Mm}{\mathcal{M}} 
\nc{\mm}{\mathrm{m}} 
\nc{\K}{\mrm{K}} 
\nc{\W}{\mrm{W}} 
\nc{\red}{{\on{red}}}
\nc{\Voev}{{\on{Voev}}}
\nc{\Corr}{\mrm{Corr}}
\nc{\Span}{\mathbf{Corr}{}}
\nc{\Gap}{\mrm{Gap}}
\nc{\Corrfr}{\Corr^{\fr}}
\nc{\Corrvfr}{\Corr^{\Vfr}}
\nc{\Spec}{\on{Spec}}
\nc{\Sm}{\on{Sm}}
\nc{\Gm}{\mathbf{G}_{\on{m}}}
\renc{\P}{\bP}
\nc{\nis}{\mathrm{nis}}
\nc{\zar}{\mathrm{zar}}
\nc{\et}{\mathrm{\acute et}}
\nc{\all}{\mathrm{all}}
\nc{\fold}{\mathrm{fold}}
\nc{\Fun}{\mathrm{Fun}}
\nc{\Ho}{\mathrm{Ho}}
\nc{\Segal}{\mathrm{Segal}}
\nc{\Mon}{\mrm{Mon}{}}
\nc{\Ab}{\mrm{Ab}}
\nc{\Sh}{\on{Sh}}
\nc{\M}{\mrm{M}}
\nc{\Lhtp}{L_{\A^1}}
\nc{\Lmot}{L_{\mrm{mot}}}
\nc{\mot}{\mrm{mot}}
\nc{\SH}{\mbf{SH}}
\nc{\RR}{\mbf{R}}
\nc{\CC}{\mbf{C}}
\nc{\Mod}{\mbf{Mod}}
\nc{\QCoh}{\mbf{QCoh}}
\nc{\MonUnit}{\mbf{1}}
\nc{\1}{\mbf{1}}
\nc{\tr}{\on{tr}}
\nc{\cotr}{\mrm{cotr}}
\nc{\vop}{\mrm{vop}}
\nc{\fr}{{\on{fr}}}
\nc{\Ar}{\mrm{Ar}}
\nc{\Vfr}{\on{Vfr}}
\nc{\frdiff}{{\on{frdiff}}}
\nc{\frGys}{\on{frGys}}
\nc{\SHfr}{\SH^{\fr}}
\nc{\SHfrdiff}{\SH^{\frdiff}}
\nc{\SHfrGys}{\SH^{\frGys}}
\nc{\InftyCat}{(\mathrm{\infty,1)\textnormal{-}Cat}}
\nc{\TriCat}{\mathrm{TriCat}}
\nc{\Cat}{\mathrm{1\textnormal{-}Cat}}
\nc{\Th}{\on{Th}}
\def\G{\bG}
\nc{\CMon}{\mrm{CMon}{}}
\nc{\MGL}{\mrm{MGL}}
\nc{\Seg}{\mrm{Seg}{}}
\nc{\GW}{\mrm{GW}{}}
\nc{\Tw}{\mrm{Tw}}
\nc{\sslash}{/\mkern-6mu/}
\nc{\PrL}{\mrm{Pr}^\mrm{L}}
\nc{\PrR}{\mrm{Pr}^\mrm{R}}
\nc{\pr}{\mrm{pr}}
\let\phi\varphi

\nc\efr{\mrm{efr}}
\nc\nfr{\mrm{nfr}}
\nc\dfr{\mrm{fr}}
\nc\tfr{\mrm{tfr}}
\nc\Vect{\mrm{Vect}}
\nc\sVect{\mrm{sVect}}
\nc{\fix}{\mrm{fix}}
\nc{\ho}{\mrm{h}}
\nc\Mfd{\mrm{Mfd}}
\nc{\PSh}{\mrm{PSh}}
\nc{\hzmw}{H \tilde\Z{}}
\nc{\Cor}{\mrm{Cor}{}}
\nc{\cormw}{\mrm{\widetilde{Cor}}{}}
\nc{\Chw}{\mrm{\widetilde{CH}}{}}
\nc{\Ex}{\mrm{Ex}}
\nc{\BM}{\mrm{BM}}
\let\setminus-
\nc{\Pic}{\mrm{Pic}}
\nc{\pur}{\mathfrak p}

\nc{\inftyCat}{\term{$\infty$-category}}
\nc{\inftyCats}{\term{$\infty$-categories}}

\nc{\inftyOneCat}{\term{$(\infty,1)$-category}}
\nc{\inftyOneCats}{\term{$(\infty,1)$-categories}}

\nc{\inftyGrpd}{\term{$\infty$-groupoid}}
\nc{\inftyGrpds}{\term{$\infty$-groupoids}}

\nc{\inftyTop}{\term{$\infty$-topos}}
\nc{\inftyTops}{\term{$\infty$-toposes}}

\nc{\inftyTwoCat}{\term{$(\infty,2)$-category}}
\nc{\inftyTwoCats}{\term{$(\infty,2)$-categories}}

\title{Framed transfers and motivic fundamental classes}

\author[E. Elmanto]{Elden Elmanto}
\address{Department of Mathematics\\
Harvard University\\
1 Oxford St.\\
Cambridge, MA 02138\\
USA}
\email{\href{mailto:elmanto@math.harvard.edu}{elmanto@math.harvard.edu}}
\urladdr{\url{https://www.eldenelmanto.com/}}
\thanks{E.E.\ was partially supported by NSF grant DMS-1508040}

\author[M. Hoyois]{Marc Hoyois}
\address{Fakultät für Mathematik\\
Universität Regensburg\\
Universitätsstr. 31\\
93040 Regensburg\\
Germany}
\email{\href{mailto:marc.hoyois@ur.de}{marc.hoyois@ur.de}}
\urladdr{\url{http://www.mathematik.ur.de/hoyois/}}
\thanks{M.H.\ was partially supported by NSF grant DMS-1761718}

\author[A. A. Khan]{Adeel A. Khan}
\address{Fakultät für Mathematik\\
Universität Regensburg\\
Universitätsstr. 31\\
93040 Regensburg\\
Germany}
\email{\href{mailto:adeel.khan@mathematik.uni-regensburg.de}{adeel.khan@mathematik.uni-regensburg.de}}
\urladdr{\url{https://www.preschema.com}}

\author[V. Sosnilo]{Vladimir Sosnilo}
\address{Laboratory ``Modern Algebra and Applications''\\
St. Petersburg State University\\
14th line, 29B\\
199178 Saint Petersburg\\
Russia}
\email{\href{mailto:vsosnilo@gmail.com}{vsosnilo@gmail.com}}
\thanks{V.S.\ was supported by the grant of the Government of the Russian Federation
for the state support of scientific research carried out under the supervision of leading scientists,
agreement 14.W03.31.0030 dated 15.02.2018}

\author[M. Yakerson]{Maria Yakerson}
\address{Fakultät für Mathematik\\
Universität Regensburg\\
Universitätsstr. 31\\
93040 Regensburg\\
Germany}
\email{\href{mailto:maria.yakerson@ur.de}{maria.yakerson@ur.de}}
\urladdr{\url{https://www.muramatik.com}}
\thanks{M.Y.\ was supported by SFB/TR 45 ``Periods, moduli spaces and arithmetic of algebraic varieties''}

\date{\today}

\begin{document}

\begin{abstract}
We relate the recognition principle for infinite $\P^1$-loop spaces to the theory of motivic fundamental classes of Déglise, Jin, and Khan. 

We first compare two kinds of transfers that are naturally defined on cohomology theories represented by motivic spectra: the framed transfers given by the recognition principle, which arise from Voevodsky's computation of the Nisnevish sheaf associated with $\A^n/(\A^n-0)$, and the Gysin transfers defined via Verdier's deformation to the normal cone.

We then introduce the category of finite $R$-correspondences for $R$ a motivic ring spectrum,
generalizing Voevodsky's category of finite correspondences and Calmès and Fasel's category of finite Milnor--Witt correspondences. Using the formalism of fundamental classes, we show that the natural functor from the category of framed correspondences to the category of $R$-module spectra factors through the category of finite $R$-correspondences.
\end{abstract}

\maketitle

\parskip 0pt
\tableofcontents

\parskip 0.2cm


\newpage
\section{Introduction}

This paper connects two recent developments in our understanding of certain cohomology theories for schemes, namely those that are represented in the Morel–Voevodsky category of motivic spectra \cite{Morel:2003}. On the one hand, the work of Levine \cite{LevineVFC} and Déglise, Jin, and Khan \cite{DJK} develops the theory of \emph{fundamental classes} in the setting of motivic homotopy theory. This results in a vast generalization of Fulton's operations in Chow groups \cite{Fulton} to these cohomology theories. On the other hand, the work of Garkusha, Panin, Ananyevskiy, and Neshitov \cite{garkusha2014framed,hitr,agp,gnp}, building on some insights of Voevodsky \cite{voevodsky2001notes}, develops a theory of \emph{framed motives}. One achievement of their work is to give explicit models for motivic suspension spectra of smooth schemes. 

Recall that if $E\in\SH(S)$ is a motivic spectrum over a scheme $S$, there is an associated bigraded cohomology theory on smooth $S$-schemes:
\[
E^{*,*}(-)\colon \Sm_S^\op \longrightarrow \Ab.
\]
Both the theory of fundamental classes and that of framed motives imply the existence of certain transfers, called \emph{framed transfers}, in such a cohomology theory. These transfers can be encoded by an extension of $E^{*,*}(-)$ to the category $\h\Span^\fr(\Sm_S)$ of framed correspondences:
\[
\begin{tikzcd}
	\Sm_S^\op \ar{r}{E^{*,*}(-)} \ar{d} & \Ab\rlap.\\
	\h\Span^\fr(\Sm_S)^\op \ar[dashed]{ur} &
\end{tikzcd}
\]
In the first part of this paper, we show that the framed transfers produced by both theories agree.
This is nontrivial as their respective constructions are based on different geometric ideas.
 In the second part of this paper, we introduce the category $\h\Span^R(\Sm_S)$ of \emph{finite $R$-correspondences} for $R$ a motivic ring spectrum, and we construct a further interesting extension 
\[
\begin{tikzcd}[column sep=6em]
	\Sm_S^\op \ar{r}{E^{*,*}(-)} \ar{d} & \Ab\\
	\h\Span^\fr(\Sm_S)^\op \ar{ur} \ar{d} & \\
	\h\Span^R(\Sm_S)^\op \ar[dashed]{uur} &
\end{tikzcd}
\]
when $E$ is a module over $R$. 
The category $\h\Span^R(\Sm_S)$ recovers Voevodsky's category of finite correspondences when $R$ is the motivic Eilenberg–Mac Lane spectrum $H\Z$, and it recovers Calmès and Fasel's category of finite Milnor–Witt correspondences when $R=\hzmw$. Thus, our construction unifies those of Voevodsky and of Calmès–Fasel, as well as their relationship with the category of framed correspondences.

\ssec{Comparison of transfers}\label{ssec:intro-comparison}
In \cite{EHKSY1}, we introduced the $\infty$-groupoid $\Corr^\fr_S(X,Y)$ of \emph{framed correspondences} between smooth $S$-schemes $X$ to $Y$: such a correspondence is a span
\begin{equation*}
  \begin{tikzcd}
     & Z \ar[swap]{ld}{f}\ar{rd}{g} & \\
    X &   & Y
  \end{tikzcd}
\end{equation*}
where $f$ is finite syntomic, together with an equivalence $\tau\colon \sL_f\simeq 0$ in $K(Z)$. Here, $\sL_f$ is the cotangent complex of $f$ and $K(Z)$ is the K-theory space of $Z$. As $X$ and $Y$ vary, these $\infty$-groupoids form the mapping spaces of an $\infty$-category $\Span^{\fr}(\Sm_S)$. 

Framed correspondences encode an essential functoriality of cohomology theories represented by motivic spectra: if $E\in \SH(S)$ and $\alpha=(Z,f,g,\tau)\in \Corr^\fr_S(X,Y)$, there is an induced map
\[
\alpha^*\colon E(Y) \to E(X)
\]
in cohomology; here $E(X)=\Maps_{\SH(S)}(\Sigma^\infty_\T X_+,E)$ is the $E$-cohomology \emph{space} of $X$. In fact, there are several different ways to construct $\alpha^*$ that are not obviously equivalent:
\begin{enumerate}
	\item\emph{Via fundamental classes.} For any finite syntomic map $f\colon Z\to X$ between $S$-schemes, its fundamental class induces a Gysin transfer $f_!\colon E(Z,\sL_f) \to E(X)$ in twisted cohomology \cite{LevineVFC,DJK}.
	Hence, given the framed correspondence $\alpha$, we can define
	\[
	\alpha^*\colon E(Y) \xrightarrow{g^*} E(Z) \stackrel\tau\simeq E(Z,\sL_f) \xrightarrow{f_!} E(X).
	\]
	\item \emph{Via Voevodsky's Lemma.} 
	Voevodsky introduced the set $\Corr^{\efr}_S(X,Y)$ of equationally framed correspondences from $X$ to $Y$~\cite{voevodsky2001notes} and constructed a canonical map
	\begin{equation*}
		\Corr^{\efr}_S(X,Y) \to \Maps_{\SH(S)}(\Sigma^\infty_\T X_+,\Sigma^\infty_\T Y_+).
	\end{equation*}
One of the key results in \cite{EHKSY1} is that the presheaves $\Corr^\efr_S(-,Y) $ and $ \Corr^\fr_S(-,Y)$ are motivically equivalent. This implies that Voevodsky's map factors through the $\infty$-groupoid $\Corr^\fr_S(X,Y)$. In particular, $\alpha$ induces a map $\Sigma^\infty_\T X_+\to \Sigma^\infty_\T Y_+$ in $\SH(S)$, whence a map $\alpha^*\colon E(Y) \to E(X)$ in cohomology.
	\item \emph{Via framed motivic spectra.} In \cite{EHKSY1}, we constructed the $\infty$-category $\SH^\fr(S)$ of framed motivic spectra over $S$, in which the functoriality with respect to framed correspondences is hard-coded. In particular, $\alpha$ induces a morphism \[\Sigma^\infty_{\T,\fr}(\alpha)\colon \Sigma^\infty_{\T,\fr}X \to \Sigma^\infty_{\T,\fr}Y\] in $\SH^\fr(S)$.
	The reconstruction theorem of \emph{loc.\ cit.} (generalized to arbitrary schemes in \cite{Hoyois:2018aa}) gives an equivalence of $\infty$-categories $\SH^\fr(S)\simeq \SH(S)$, under which $\Sigma^\infty_{\T,\fr}(\alpha)$ corresponds to a morphism $\Sigma^\infty_\T X_+\to \Sigma^\infty_\T Y_+$ as in (2).
\end{enumerate}
In \secref{sec:transfers}, we show that these three constructions agree.

We note that each construction has its own useful features.
Construction (1) connects framed correspondences with the powerful formalism of six operations. As we explained in the introduction to \cite{EHKSY1}, it was this hypothetical connection that led us to the correct formulation of the recognition principle for infinite $\P^1$-loop spaces.
Construction (2) is helpful to perform explicit computations. For example, Bachmann and Yakerson employ the Voevodsky transfer to show that for a strictly homotopy invariant Nisnevich sheaf of abelian groups $M$ on $\Sm_k$ the double contraction $M_{-2}$ has an infinite $\G_m$-delooping (at least when $\chr k = 0$) \cite{BY}.
Construction (3) has the advantage that it is coherently compatible with the \emph{composition} of framed correspondences, i.e., it gives a functor
\[
\Span^\fr(\Sm_S) \to \SH(S).
\]
Our comparison theorems can therefore also be viewed as coherence theorems for the first two types of transfers.

\ssec{Finite correspondences for motivic ring spectra}\label{ssec:intro-CorrE}
In~\secref{sec:fin-corr}, we introduce categories of finite correspondences that encode the functoriality of $R$-cohomology for a given motivic ring spectrum $R\in\SH(S)$. We define for $X,Y\in\Sm_S$ an $\infty$-groupoid $\Corr^R_S(X,Y)$ of \emph{finite $R$-correspondences} such that $\Corr^R_S(X,S) \simeq R(X)$. We expect that these are the mapping spaces of an $\infty$-category $\Span^R(\Sm_S)$ with the following properties:
\begin{enumerate}
	\item There is a functor $M^R\colon \Span^R(\Sm_S) \to \Mod_R(\SH(S))$ sending $X$ to $R\otimes \Sigma^\infty_\T X_+$.
	\item There is a functor $\Phi^R\colon \Span^\fr(\Sm_S) \to \Span^R(\Sm_S)$ sending $X$ to $X$.
	\item The following square of $\infty$-categories commutes:
	\[
	\begin{tikzcd}
		\Span^\fr(\Sm_S) \ar{r} \ar{d}{\Phi^R} & \SH(S) \ar{d}{R\otimes -} \\
		\Span^R(\Sm_S) \ar{r}{M^R} & \Mod_R(\SH(S)). 
	\end{tikzcd}
	\]
\end{enumerate}
In this paper we restrict ourselves to constructing the homotopy category $\ho\Span^R(\Sm_S)$, and we establish the above properties at the level of homotopy categories (enriched in the homotopy category of spaces). 
The functor $M^R$ exists essentially by design, and the functor $\Phi^R$ is defined using the formalism of fundamental classes. Property (3) is an application of the main result of \secref{sec:transfers}.

We then consider the cases $R=H\Z$ and $R=\hzmw$ for $S$ essentially smooth over a Dedekind domain and a field, respectively. In these two cases the mapping spaces $\Corr^R_S(X,Y)$ are discrete, so the $\infty$-category $\Span^R(\Sm_S)$ is defined and is a $1$-category. Moreover, we prove the following comparison results:
\begin{enumerate}\setcounter{enumi}{3}
	\item $\Span^{H\Z}(\Sm_S)$ is equivalent to Voevodsky's category of finite correspondences, and the functor $\Phi^{H\Z}$ is the functor $\cyc$ constructed in \cite[\sectsign 5.3]{EHKSY1}.
	\item $\Span^{\hzmw}(\Sm_S)$ is equivalent to the category of finite Milnor–Witt correspondences, constructed by Calmès and Fasel \cite{Calmes:2014ab}, and the functor $\Phi^{\hzmw}$ refines the one defined by Déglise and Fasel in \cite{DegliseFasel}.
\end{enumerate}

If $k$ is a field and $R\in \SH(k)$ is an $\mathrm{MSL}$-algebra, the category $\ho\Span^R(\Sm_k)$ is equivalent to that constructed by Druzhinin and Kolderup in \cite{DruzhininKolderup}. For $R=\mathrm{KGL}$ (resp.\ $R=\mathrm{KO}$ if $\chr k\neq 2$), it receives a functor from Walker's category of finite $K_0$-correspondences \cite{SuslinGrayson} (resp.\ from Druzhinin's category of finite $\GW$-correspondences \cite{DruzhininGW}).
However, a novel feature of our category is that it is enriched in the homotopy category of spaces, hinting that it is the homotopy category of a more fundamental $\infty$-category. Its mapping spaces are discrete if and only if $R$ is $0$-truncated in the effective homotopy $t$-structure, a condition which implies that $R$ is an $\hzmw$-algebra. The Calmès–Fasel category $\Span^{\hzmw}(\Sm_k)$ is thus in a precise sense the most general \emph{1-category} of finite correspondences.

	Assuming that the $\infty$-category $\Span^R(\Sm_S)$ has been constructed, one can consider the $\infty$-category $\DM^R(S)$ of $\T$-spectra in $\A^1$-invariant Nisnevich-local presheaves on $\Span^R(\Sm_S)$.
	When $S$ is the spectrum of a field of characteristic zero and $R=H\Z$ or $R=\hzmw$, it is well known that $\DM^R(S)\simeq\Mod_R(\SH(S))$ \cite{Rondigs:2008,ElmantoKolderup} and that the ``cancellation theorem'' holds for $\Span^R(\Sm_S)$ \cite{voevodsky-cancel,FaselOstvaer}.
	 We will not attempt here to generalize these results. However, we note that the conjectural properties listed above imply that $\SH(S)$ is always a retract of $\DM^\1(S)$.
	 
\ssec{Conventions and notation}
Our terminology and notation follows \cite{EHKSY1}. In particular:
\begin{itemize}
	\item $\Spc$ is the $\infty$-category of spaces/$\infty$-groupoids, $\Spt$ that of spectra;
	\item $\Maps(X,Y)$ is the space of maps from $X$ to $Y$ in an $\infty$-category;
	\item if $C$ is an $\infty$-category, we denote by $\ho C$ its homotopy category;
	\item $\Perf(X)$ is the $\infty$-category of perfect complexes over $X$;
	\item $\SH(S)$ is the stable motivic homotopy $\infty$-category over $S$;
	\item $\DM(S)$ is Voevodsky's $\infty$-category of motives over $S$;
	\item $\G$, $\T$, and $\P$ denote the pointed presheaves $(\G_m,1)$, $\A^1/\G_m$, and $(\P^1,\infty)$; $\Sigma_\G$, $\Sigma_\T$, and $\Sigma_\P$ are the corresponding suspension functors, and $\Omega_\G$, $\Omega_\T$, and $\Omega_\P$ their right adjoints.
	\item If $C$ is an $\infty$-category and $M, N$ are two collections of morphisms in $C$ that are stable under composition and pullback along one another, we write $\Span(C, M, N)$ for the $\infty$-category of spans with backward maps in $M$ and forward maps in $N$; see \cite[\sectsign 5]{BarwickMackey} for details on the construction of this $\infty$-category.
\end{itemize}

\ssec{Acknowledgments}

We would like to thank the Mittag-Leffler Institute and the organizers of the program ``Algebro-Geometric and Homotopical Methods'', which hosted the authors while part of this work was done.
Elmanto and Yakerson also thank the University of Copenhagen and the Centre for Symmetry and Deformation for support and hospitality during an enjoyable research visit.
Hoyois would like to thank Aravind Asok, Tom Bachmann, and Marc Levine for useful discussions about Chow–Witt groups.
Yakerson would like to thank sincerely Jean Fasel and Marc Levine for teaching her about Milnor--Witt correspondences.


\section{Preliminaries}

In this section, we review some aspects of the formalism of six functors in stable motivic homotopy theory \cite{Ayoub,CD}.

In \ssecref{ssec:diff/review}, we discuss various (co)homology theories associated with a motivic spectrum and their basic properties. In \ssecref{ssec:gysin}, we review the formalism of fundamental classes for local complete intersection morphisms.


\ssec{Cohomology theories}
\label{ssec:diff/review}

\sssec{} 
For every morphism of schemes $f\colon X \rightarrow Y$ we have an adjunction
\[
f^*: \SH(Y) \rightleftarrows \SH(X): f_*.
\]
If $f$ is smooth, there is a further left adjoint
\[
f_\sharp: \SH(X) \rightleftarrows \SH(Y): f^*.
\]
If $f$ is locally of finite type\footnote{The careful reader will replace ``locally of finite type'' by ``separated of finite type'', since the current literature only contains the construction of $f_!$ in the latter case; it can be constructed in the general case using Zariski descent.}, we also have an adjunction
\[
f_!: \SH(X) \rightleftarrows \SH(Y):f^!,
\]
such that $f_!\simeq f_*$ if $f$ is proper.

The basic properties of these functors are summarized by the existence of a functor
\[
\Span(\Sch,\mathrm{all},\mathrm{lft}) \to \InftyCat, \quad S\mapsto \SH(S), \quad (U \xleftarrow f T \xrightarrow p S) \mapsto p_!f^*,
\]
where ``$\mathrm{lft}$'' is the class of morphisms locally of finite type (see \cite[Chapter 2, \sectsign 5.2]{KhanThesis} or \cite[\sectsign 6.2]{hoyois-sixops}). We will often use this implicitly when discussing the functoriality of certain constructions.

\sssec{Thom transformations}
Let $S$ be a scheme, $\sE$ a locally free $\sO_S$-module of finite rank, and $V=\Spec(\Sym(\sE))$ the associated vector bundle.
If $p\colon V\rightarrow S$ is the structure morphism and $s\colon S \rightarrow V$ the zero section, then the adjoint functors
\[
\Sigma^{\sE}=p_{\sharp}s_*: \SH(S) \rightleftarrows \SH(S): s^!p^*=\Sigma^{-\sE}
\]
are $\SH(S)$-linear equivalences of $\infty$-categories, called \emph{Thom transformations}. In particular, $\Sigma^\sE\simeq \Sigma^\sE\1_S\otimes(-)$, and the object $\Sigma^\sE\1_S\in \SH(S)$ is invertible with inverse $\Sigma^{-\sE}\1_S$.

The Thom transformations $\Sigma^\sE$ are defined more generally for $\sE$ a perfect complex of $\sO_S$-modules, and they assemble into a morphism of grouplike $\Einfty$-spaces
\begin{equation} \label{k-to-sh}
K(S) \rightarrow \Pic(\SH(S)),\quad  \xi \mapsto \Sigma^\xi\1_S,
\end{equation}
natural in $S$, called the \emph{motivic J-homomorphism} (see \cite[\sectsign 16.2]{norms}).
In particular, for every cofiber sequence
\[
\mathcal E' \rightarrow \mathcal E \rightarrow \mathcal E''
\]
in $\Perf(S)$, we have canonical equivalences
\begin{equation} \label{exact-thom}
\Sigma^\mathcal E \simeq \Sigma^{\mathcal E'} \Sigma^{\mathcal E''} \simeq \Sigma^{\mathcal E''} \Sigma^{\mathcal E'}.
\end{equation}

\sssec{Purity equivalences}
For $f\colon X \rightarrow S$ a smooth morphism with sheaf of relative differentials $\Omega_f$, we have canonical equivalences
\begin{equation} \label{purity1}
f^!\simeq \Sigma^{\Omega_f}f^*\quad\text{and}\quad f_! \simeq f_{\sharp}\Sigma^{-\Omega_f}.
\end{equation}
Suppose that $s\colon Z\hookrightarrow X$ is a closed immersion such that the composite $g=f\circ s$ is smooth. Combining \eqref{exact-thom} and \eqref{purity1}, we obtain equivalences
\begin{equation}\label{purity2}
	s^!f^* \simeq \Sigma^{-\sN_s} g^*\quad\text{and}\quad f_\sharp s_* \simeq g_\sharp\Sigma^{\sN_s},
\end{equation}
where $\sN_s$ is the conormal sheaf of $s$. The equivalences \eqref{purity1} and \eqref{purity2} are called the \emph{purity equivalences}. Note that we have equivalences of perfect complexes $\Omega_f\simeq \sL_f$ and $\sN_s[1]\simeq \sL_s$, so that we can write $\Sigma^{\Omega_f}\simeq \Sigma^{\sL_f}$ and $\Sigma^{-\sN_s}\simeq \Sigma^{\sL_s}$.

\sssec{Twisted cohomology}\label{sssec:theories}
A motivic spectrum $E\in \SH(S)$ gives rise to various (co)homology theories for $S$-schemes, which can be twisted by $K$-theory classes.
Let $p\colon X\to S$ be a morphism\footnote{Whenever the functors $p_!$ or $p^!$ are used, it is implicitly assumed that $p$ is locally of finite type.} and let $\xi\in K(X)$. We will consider the following mapping spaces:
\begin{enumerate}
	\item The \emph{$\xi$-twisted cohomology} of $X$ with coefficients in $E$ is
	\[
	E(X,\xi) = \Maps(\1_S, p_*\Sigma^{\xi} p^* E).
	\]
	\item The \emph{$\xi$-twisted Borel–Moore homology} of $X$ with coefficients in $E$ is
	\[
	E^\BM(X/S,\xi) = \Maps(\1_S, p_*\Sigma^{-\xi} p^! E).
	\]
\end{enumerate}
We omit the second parameter when $\xi=0$. Moreover, it is understood that an element $\xi\in K(X)$ is allowed to twist the cohomology of any $X$-scheme: if $f\colon X'\to X$ is a morphism, we will often write $E(X',\xi)$ instead of $E(X',f^*\xi)$, and similarly for Borel–Moore homology.
 There are also twisted versions of compactly supported cohomology and of homology (see \cite[Definition 2.2.1]{DJK}), but we shall not use these theories in this paper.

\begin{rem}\label{rem:E}
	In what follows, we often fix a motivic spectrum $E\in\SH(S)$ and talk about $E$-cohomology spaces in the interest of readability. However, $E$-cohomology spaces can generally be replaced by the corresponding endofunctors of $\SH(S)$. In particular, the naturality in $E$ of all constructions and statements will be implicit.
\end{rem}

\sssec{Twisted motives}
\label{sssec:homological-motives}
One can also define various \emph{twised motives} in $\SH(S)$: if $p\colon X\to S$ is a morphism and $\xi\in K(X)$, we let
	\[
	M_S(X,\xi) = p_!\Sigma^\xi p^!\1_S\quad\text{and}\quad M_S^\BM(X,\xi)=p_!\Sigma^{-\xi} p^*\1_S.
	\]
	
	For every $E\in\SH(S)$, we have
	\[
	E^\BM(X/S,-\xi) \simeq \Maps(M_S^\BM(X,\xi), E)
	\]
	by adjunction. 
	The relationship between $M_S(X,\xi)$ and cohomology is more subtle.
	There is a canonical map $p^*E\to \Hom(p^!\1_S, p^!E)$ adjoint to the composite
	\[
	p_!(p^*E\otimes p^!\1_S) \simeq E \otimes p_!p^!\1_S \xrightarrow{\text{counit}} E\otimes\1_S\simeq E.
	\]
	Applying $\Maps(\Sigma^{\xi}\1_X,-)$, we obtain a canonical map
	\[
	E(X,-\xi) \to \Maps(M_S(X,\xi), E);
	\]
	it is an equivalence when $X$ is smooth over $S$ by purity \eqref{purity1}, whence when $X$ is cdh-locally smooth since motivic spectra satisfy cdh descent \cite[Proposition 3.7]{Cisinski}. However, it is not known to be an equivalence in general.

\sssec{Functoriality}\label{sssec:functoriality}
The cohomology space $E(X,\xi)$ is contravariant in the pair $(X,\xi)$. More precisely, if $(\Sch_S)_{/K}\to \Sch_S$ denotes the Cartesian fibration classified by $K\colon \Sch_S^\op \to \Spc$, then $(X,\xi)\mapsto E(X,\xi)$ is a contravariant functor on $(\Sch_S)_{/K}$. In particular, for every $S$-morphism $f\colon Y\to X$, there is a pullback map
\[
f^*\colon E(X,\xi) \to E(Y,f^*\xi)
\]
induced by the unit transformation $\id\to f_*f^*$. 

On the other hand, Borel–Moore homology $E^\BM(X/S,\xi)$ is covariant in $(X,\xi)$ for proper maps and contravariant in $(X,\xi)$ for étale maps. This bivariance can be expressed coherently using the $\infty$-category of correspondences $\Span((\Sch_S)_{/K},\mathrm{prop},\et)$. In addition, Borel–Moore homology is contravariantly functorial in the base $S$. In particular, for a morphism $f\colon S'\to S$, there is a base change map
\[
f^*\colon E^\BM(X/S,\xi) \to E^\BM(X\times_SS'/S', \pi_1^*\xi)
\]
induced by the exchange transformations $\Ex^*_*$ and $\Ex^{*!}$.

\sssec{Cohomology with support} 
\label{sssec:E_Z(X)}

Let $X$ be an $S$-scheme, $i\colon Y\hookrightarrow X$ an immersion, and $\xi\in K(Y)$. The \emph{$\xi$-twisted cohomology of $X$ with support in $Y$} is
\[
E_Y(X,\xi) = \Maps(\1_S, p_* i_!\Sigma^{\xi} i^!p^*E).
\]

Given a Cartesian square
\[
\begin{tikzcd}
	Y' \ar[hook]{r}{j} \ar{d}[swap]{g} & X' \ar{d}{f} \\
	Y \ar[hook]{r}{i} & X,
\end{tikzcd}
\]
the unit transformation $\id\to g_*g^*$ and the exchange transformations $\Ex^{*!}\colon g^* i^! \to j^! f^*$ and $\Ex_{!*}\colon i_!g_* \to f_*j_!$ define a transformation
\[
i_!\Sigma^{\xi}i^! \rightarrow f_*j_!\Sigma^{g^*\xi}j^!f^*,
\]
which induces a pullback in cohomology with support
\begin{equation} \label{eqn:pull-w-support}
f^*\colon E_Y(X,\xi) \to E_{Y'}(X',g^*\xi).
\end{equation}

If $k\colon V\hookrightarrow Y$ is another immersion, we also have a ``forgetful'' map
\begin{equation}\label{eqn:forget-support}
E_V(X,k^*\xi) \to E_Y(X,\xi)
\end{equation}
induced by the counit transformation $k_!k^!\to \id$.

If $Y$ is closed in $X$ and both are smooth over a common base, we have a purity equivalence
\begin{equation}\label{eqn:purity-support}
E_Y(X,\xi) \simeq E(Y,\xi-\sN_i)
\end{equation}
by~\eqref{purity2}.

\sssec{Localization}
Suppose that we have a diagram in $\Sch_S$
\begin{equation*}
 Z \stackrel{i}{\hookrightarrow} X \stackrel{j}{\hookleftarrow} U
\end{equation*}
where $i$ is a closed immersion and $j$ is the complementary open immersion, and let $\xi\in K(X)$.
Then the localization sequence
\[
i_!i^! \rightarrow \id \rightarrow j_*j^*
\]
gives the fiber sequence
\begin{equation}\label{eqn:E_Z(X)}
E_Z(X,i^*\xi) \rightarrow E(X,\xi) \rightarrow E(U,j^*\xi).
\end{equation}
Dually, the localization sequence
\[
j_!j^! \rightarrow \id \rightarrow i_*i^*
\]
gives the fiber sequence
\begin{equation}\label{eqn:E(X/Z)}
E_U(X,j^*\xi) \rightarrow E(X,\xi) \rightarrow E(Z,i^*\xi).
\end{equation}

In Borel–Moore homology, we similarly obtain the fiber sequence
\[
E^\BM(Z/S,i^*\xi) \rightarrow E^\BM(X/S,\xi) \rightarrow E^\BM(U/S,j^*\xi).
\]

\sssec{Descent properties}
\label{sssec: descent properties}
Recall that the functor
\[
\Sch_S^\op \to \Fun(\SH(S),\SH(S)),\quad (p\colon X\to S)\mapsto p_*p^*,
\]
is an $\A^1$-invariant cdh-sheaf on $\Sch_S$ \cite[Proposition 3.7]{Cisinski}. 
Consequently, cohomology is an $\A^1$-invariant cdh-sheaf and Borel–Moore homology is a Nisnevich sheaf. In particular, if $f\colon X \rightarrow Y$ is an $\A^1$-cdh-equivalence (i.e., $f$ induces an equivalence between the associated $\A^1$-invariant cdh sheaves) and $\xi\in K(Y)$, then the induced map
\[
f^*\colon E(Y,\xi) \to E(X,\xi)
\]
is an equivalence. 

In fact, we have the following more precise excision properties.
Let $Y\subset X$ be a subscheme and $f\colon X'\to X$ a morphism such that $f^{-1}(Y)\simeq Y$. For any $\xi\in K(Y)$, the pullback
\[
f^*\colon E_Y(X,\xi) \to E_{Y}(X',\xi)
\]
is an equivalence under either of the following conditions:
\begin{itemize}
	\item $f$ is smooth and $Y$ is closed (Nisnevich excision);
	\item $f$ is proper and $Y$ is open (excision for abstract blowups).
\end{itemize}
This follows directly from the definition of $f^*$ given in \sssecref{sssec:E_Z(X)}.

\sssec{Products}\label{sssec:products}

Suppose that $E\in\SH(S)$ is equipped with a multiplication $\mu\colon E\otimes E\to E$. This induces various products in cohomology and in Borel–Moore homology:
\begin{enumerate}
	\item For any $S$-scheme $X$, subschemes $Z,Z'\subset X$, and $K$-theory classes $\xi\in K(Z)$ and $\xi'\in K(Z')$, we have the usual \emph{cup product}
\[
\mu\colon E_Z(X,\xi) \times E_{Z'}(X,\xi') \to E_{Z\cap Z'}(X,\xi+\xi'),\quad (x,y)\mapsto x\cup y.
\]
	\item For any $S$-scheme $X$, subschemes $T\subset Z\subset X$, and $K$-theory classes $\xi\in K(Z)$ and $\zeta\in K(T)$, we have the \emph{refined cup product}
\[
\bar\mu\colon E_Z(X,\xi) \times E_T(Z,\zeta) \to E_T(X,\xi+\zeta).
\]
We refer to \cite[1.2.8]{deglise2011orientation} for the definition.
This refines the cup product from (1) as follows: there is a commutative square
\[
\begin{tikzcd}
	E_Z(X,\xi) \times E_{Z'}(X,\xi') \ar{r}{\id\times i^*} \ar{d}[swap]{i'^*\times\id} \ar{dr}[description]{\mu} & E_Z(X,\xi) \times E_{Z\cap Z'}(Z,\xi') \ar{d}{\bar\mu} \\
	E_{Z\cap Z'}(Z',\xi) \times E_{Z'}(X,\xi') \ar{r}[swap]{\bar\mu} & E_{Z\cap Z'}(X,\xi+\xi'),
\end{tikzcd}
\]
where $i\colon Z\hook X$ and $i'\colon Z'\hook X$ are the inclusions.
	\item Suppose $Z\to Y\to X$ are $S$-morphisms locally of finite type and let $\xi\in K(Z)$ and $\zeta\in K(Y)$. Then we have the \emph{composition product}
\[
\mu^\BM\colon E^\BM(Z/Y,\xi) \times E^\BM(Y/X,\zeta) \to E^\BM(Z/X,\xi+\zeta),\quad (z,y)\mapsto z\cdot y.
\]
We refer to \cite[1.2.8]{deglise2017bivariant} for the definition.
\end{enumerate}
Of course, the cup product and the composition product are associative or unital (up to homotopy) if the multiplication on $E$ is.

\sssec{Borel–Moore homology as cohomology with support}
\label{sssec:bm-support}
Let $f\colon Z\to S$ be a morphism locally of finite type. We say that $f$ is \emph{smoothable} if there exists a factorization
\[
\begin{tikzcd}
	Z \ar[hook]{r}{i} \ar{dr}[swap]{f} & X \ar{d}{p} \\
	& S
\end{tikzcd}
\]
where $i$ is a closed immersion and $p$ is smooth. For example, if $S$ has the resolution property (i.e., every finitely generated quasi-coherent sheaf is a quotient of a locally free sheaf of finite rank), then every quasi-projective morphism $f\colon Z\to S$ is smoothable.

In the above situation, if $E\in \SH(S)$ and $\xi\in K(Z)$, the purity equivalence~\eqref{purity1} induces a canonical equivalence
\begin{equation}\label{purity3}
E^\BM(Z/S,\xi) \simeq E_Z(X,\Omega_{X/S}-\xi).
\end{equation}
We record the following compatibility properties of the equivalence~\eqref{purity3}, which follow easily from the definitions. We state them without twists for simplicity.
\begin{enumerate}
	\item \emph{Base change.} The equivalence~\eqref{purity3} is contravariantly functorial in $S$.
	\item \emph{Pushforwards.} Consider a commutative diagram
	\[
	\begin{tikzcd}
		T \ar[hook]{r}{k} \ar{d}[swap]{h} & Y \ar{d}{g} \\
		Z \ar[hook]{r}{i} \ar{dr} & X \ar{d}{f} \\
		& S
	\end{tikzcd}
	\]
	where $f$ and $g$ are smooth, $h$ is proper, and $i$ and $k$ are closed immersions. Then the following diagram commutes:
	\[
	\begin{tikzcd}
		E^\BM(T/S) \ar{r}{h_*} \ar{d}[swap]{\simeq} & E^\BM(Z/S) \ar{d}{\simeq} \\
		E_T(Y,\Omega_{Y/S}) \ar{r}{g_!} & E_Z(X,\Omega_{X/S}).
	\end{tikzcd}
	\]
	Here, $h_*$ is the proper pushforward and $g_!$ is the Gysin map induced by the purity equivalence for $g$ (see \sssecref{sssec:pushforwad-coh} for the definition of $g_!$ in a more general context). As a special case, if $t\colon W\hook Z$ is a closed immersion, then the proper pushforward $t_*\colon E^\BM(W/S)\to E^\BM(Z/S)$ is identified with the forgetful map $E_W(X,\Omega_{X/S})\to E_Z(X,\Omega_{X/S})$.
	\item \emph{Products.} Suppose that $E$ is equipped with a multiplication $\mu\colon E\otimes E\to E$, and consider a commutative diagram
	\[
	\begin{tikzcd}
		T \ar[hook]{r} \ar{dr} & V \ar[hook]{r} \ar{d} & Y \ar{d}{q} \\
		& Z \ar[hook]{r} \ar{dr} & X \ar{d}{p} \\
		& & S
	\end{tikzcd}
	\]
	where the vertical maps are smooth, the horizontal maps are closed immersions, and the square is Cartesian. Then the following diagram commutes:
	\[
	\begin{tikzcd}
		E^\BM(T/Z) \times E^\BM(Z/S) \ar{r}{\mu^\BM} \ar{d}[swap]{\simeq} & E^\BM(T/S) \ar{dd}{\simeq} \\
		E_T(V,\Omega_{V/Z}) \times E_Z(X,\Omega_{X/S}) \ar{d}[swap]{\id\times q^*} & \\
		E_T(V,\Omega_{V/Z}) \times E_V(Y, q^*\Omega_{X/S}) \ar{r}{\bar\mu} & E_T(Y,\Omega_{Y/S}).
	\end{tikzcd}
	\]
\end{enumerate}


\ssec{Fundamental classes}
\label{ssec:gysin}

\sssec{}\label{sssec:purity-transformation}
We briefly recall the formalism of fundamental classes from \cite{DJK}.
Let $f\colon X\to Y$ be a smoothable
lci morphism. The \emph{fundamental class} of $f$ is a canonical element 
\[
\eta_f\in\pi_0\1^\BM(X/Y,\sL_f)=\pi_0\Maps(\Sigma^{\sL_f}\1_X,f^!\1_Y).
\]
 The associated \emph{purity transformation}
\[
\pur_f\colon \Sigma^{\sL_f} f^* \to f^!
\]
is defined as the composition
\[
\Sigma^{\sL_f} f^*(E) \simeq  \Sigma^{\sL_f} \1_X \otimes f^*(E) \xrightarrow{\eta_f\otimes \id}  f^!(\1_Y)\otimes f^*(E) \to f^!(\1_Y\otimes E) \simeq f^!(E),
\]
where the last morphism is the canonical one (see for example \cite[2.1.10]{DJK}).
The following proposition summarizes the key properties of fundamental classes:

\begin{prop}\label{prop:Gysin}\leavevmode

\noindent{\em(i)}
Let $f\colon X \rightarrow Y$ and $g\colon Y \rightarrow Z$ be lci morphisms such that $g$, $g\circ f$, and hence $f$ are smoothable. Then the following diagram commutes:
\begin{equation*}
\begin{tikzcd}
\Sigma^{\sL_{g \circ f}}(g \circ f)^* \ar{r}{\pur_{g \circ f}} \ar{d}{\simeq} & (g \circ f)^! \ar{d}{\simeq} \\
\Sigma^{\sL_f}f^*\Sigma^{\sL_g} g^* \ar{r}{\pur_f \pur_g} & f^! g^!.
\end{tikzcd}
\end{equation*}
Here, the left vertical arrow uses the equivalence $\sL_{g\circ f}\simeq f^*(\sL_g)+ \sL_f$ in $K(X)$ induced by the canonical cofiber sequence $f^*(\sL_g)\to \sL_{g\circ f}\to \sL_f$ in $\Perf(X)$.

\noindent{\em(ii)} Given a tor-independent Cartesian square
 \begin{equation*} \label{eq:tor-indep}
\begin{tikzcd}
X' \ar{d}{g} \ar{r}{v} & X \ar{d}{f} \\
Y' \ar{r}{u} & Y
\end{tikzcd}
\end{equation*}
where $f$ is lci and smoothable, the following diagrams commute:
\begin{equation*} \label{eq:*!-vfc}
\begin{tikzcd}
v^*\Sigma^{\sL_f}f^* \ar{d}{\simeq} \ar{r}{\pur_f} & v^*f^! \ar{d}{\Ex^{*!}} \\
\Sigma^{\sL_g} g^*u^* \ar{r}{\pur_g} & g^!u^*,
\end{tikzcd}
\qquad
\begin{tikzcd}
v^!\Sigma^{\sL_f}f^* \ar[leftarrow]{d}{\Ex^{*!}} \ar{r}{\pur_f} & v^!f^! \ar[leftarrow]{d}{\simeq} \\
\Sigma^{\sL_g} g^*u^! \ar{r}{\pur_g} & g^!u^!.
\end{tikzcd}
\end{equation*}
Here, the left vertical arrows use the equivalence $v^*(\sL_f) \simeq \sL_g$ in $\Perf(X')$.

\noindent{\em(iii)} If $f\colon X \rightarrow S$ is smooth, then $\pur_f\colon \Sigma^{\sL_f}f^* \rightarrow f^!$ coincides with the purity equivalence~\eqref{purity1}.
\end{prop}

\begin{proof}
For assertions (i) and (ii), see \cite[Propositions 2.5.4 and 2.5.6]{DJK}. Assertion (iii) holds by construction of $\eta_f$, see \cite[Theorem 3.3.2(1)]{DJK}.
\end{proof}

As a consequence of Proposition~\ref{prop:Gysin}(i,iii), if $f\colon X\to S$ is smooth and $i\colon Z\hook X$ is a closed immersion such that $f\circ i$ is smooth, the transformation
\[
\pur_i f^*\colon \Sigma^{\sL_i} i^* f^* \to i^! f^*
\]
coincides with the purity equivalence~\eqref{purity2}.

\sssec{Gysin maps in cohomology}
\label{sssec:pushforwad-coh}
Consider a commutative square of $S$-schemes
\[
\begin{tikzcd}
	Z \ar[hook]{r}{i} \ar{d}[swap]{g} & X \ar{d}{f} \\
	T \ar[hook]{r}{k} & Y\rlap,
\end{tikzcd}
\]
where $f$ is smoothable and lci, $i$ and $k$ are closed immersions, and $g$ is proper.\footnote{More generally, it suffices to assume that $i$ and $k$ are immersions and that the scheme-theoretic image of $Z$ in $X$ is proper over $Y$, so that $f_*i_!\simeq f_!i_!$. This makes~\eqref{eqn:forget-support} a special case of~\eqref{eq:gysin-in-coh}.} For every $\xi\in K(T)$, we have a pushforward morphism or \emph{Gysin map}
\begin{equation} \label{eq:gysin-in-coh}
f_!\colon E_{Z}(X,g^*\xi+i^*\sL_f) \to E_{T}(Y,\xi),
\end{equation}
defined by the composition
\[
f_*i_!\Sigma^{g^*\xi+i^*\sL_f}i^!f^* \xrightarrow{\pur_f}
f_*i_!\Sigma^{g^*\xi}i^!f^! \simeq k_!g_!g^!\Sigma^{\xi}k^! \xrightarrow{\text{counit}} k_!\Sigma^\xi k^!.
\]
Let us emphasize two special cases:
\begin{enumerate}
	\item If $i\colon Z\hookrightarrow X$ is a regular closed immersion and $\xi\in K(Z)$, we have the Gysin map $i_!\colon E(Z,\xi) \to E_Z(X,\xi-\sL_i)$, which generalizes the equivalence~\eqref{eqn:purity-support}.
	\item If $f\colon X\to Y$ is smoothable, lci, and proper, and if $\xi\in K(Y)$, we have the Gysin map $f_!\colon E(X,f^*\xi+\sL_f) \to E(Y,\xi)$.
\end{enumerate}

Properties (i) and (ii) of Proposition~\ref{prop:Gysin} imply obvious compatibilities of these Gysin maps with composition and pullback.

\sssec{Gysin maps in Borel–Moore homology}
\label{sssec:pullback-BM}
Let $f\colon X\to Y$ be a smoothable lci $S$-morphism. For every $\xi\in K(Y)$, there is a pullback morphism or \emph{Gysin map}
\[
f^!\colon E^\BM(Y/S,\xi) \to E^\BM(X/S,f^*\xi+\sL_f),
\]
defined by the composition
\[
\Sigma^{-\xi} \xrightarrow{\text{unit}} f_*f^*\Sigma^{-\xi} \simeq f_*\Sigma^{-f^*\xi} f^* \xrightarrow{\pur_f} f_*\Sigma^{-f^*\xi-\sL_f}f^!.
\]

Properties (i) and (ii) of Proposition~\ref{prop:Gysin} imply obvious compatibilities of these Gysin maps with composition, proper pushforward, and tor-independent base change.

\begin{rem}
We use the notation $f_!$ and $f^!$ for Gysin maps, rather than $f_*$ and $f^*$, as a visual reminder that these maps use the purity transformation $\pur_f$.
It does not indicate a particular relation to the functors $f_!$ and $f^!$.
\end{rem}

\sssec{Functoriality}\label{sssec:eta-functoriality}
If $M$ is a collection of morphisms of schemes that is closed under tor-independent base change, we let \[\Fun^{\mathrm{cart},\pitchfork}_{M}(\Delta^1,\Sch)\subset\Fun(\Delta^1,\Sch)\] be the subcategory whose objects are the morphisms in $M$ and whose morphisms are the tor-independent Cartesian squares. By Proposition~\ref{prop:Gysin}(ii), the assignment $f\mapsto \eta_f$ is a section of the Cartesian fibration classified by the functor
\[
\Fun^{\mathrm{cart},\pitchfork}_{\mathrm{slci}}(\Delta^1,\Sch)^\op \to \tau_{\leq 0}\Spc, \quad (f\colon X\to Y)\mapsto \tau_{\leq 0}\1^\BM(X/Y,\sL_f),
\]
where ``$\mathrm{slci}$'' is the collection of smoothable lci morphisms.

We expect that the construction $f\mapsto \eta_f$ can be refined to a section of
\[
\Fun^{\mathrm{cart},\pitchfork}_{\mathrm{slci}}(\Delta^1,\Sch)^\op \to \Spc, \quad (f\colon X\to Y)\mapsto \1^\BM(X/Y,\sL_f),
\]
but this is a nontrivial task because the construction of $\eta_f$ depends on a choice of factorization of $f$. For our purposes, it will suffice to know that we do have such refinements on the subcategory of regular immersions or that of smooth morphisms. In the case of regular closed immersions, the construction of the fundamental class in \cite[\sectsign 3.2]{DJK} is clearly functorial, since blowing-up commutes with tor-independent base change.
The case of regular immersions follows since an immersion factors canonically as a closed immersion followed by an open immersion.
The case of smooth morphisms can be reduced to the case of regular immersions by expressing the purity transformation for a smooth morphism in terms of the purity transformation for its diagonal, as in \cite[(2.3.4.a)]{DJK}.

The functoriality of the fundamental class $f\mapsto \eta_f$ propagates to the purity transformation $f\mapsto\pur_f$ and to the Gysin map $f\mapsto f_!$. For example, the Gysin map for regular closed immersions can be viewed as a natural transformation between the two functors
\begin{align*}
	\Fun^{\mathrm{cart},\pitchfork}_{\mathrm{reg.cl.imm}}(\Delta^1,\Sch_S)^\op &\to \Spc, \\
	(i\colon Z\hook X) &\mapsto E(Z),\\
	(i\colon Z\hook X) &\mapsto E_Z(X,-\sL_i). 
\end{align*}

\sssec{} \label{sssec:eta-composition}
We now discuss the functoriality of the commutative square of Proposition~\ref{prop:Gysin}(i).
Let \[\Fun_{M_0,M_1,M_2}^{\mathrm{cart},\pitchfork}(\Delta^2,\Sch) \subset \Fun(\Delta^2,\Sch)\] be the subcategory whose objects are triangles
\[
\begin{tikzcd}
	Z \ar{dr}[swap]{f_1} \ar{r}{f_2} & Y \ar{d}{f_0} \\ & X
\end{tikzcd}
\]
with $f_i\in M_i$ and whose morphisms are natural transformations composed of tor-independent Cartesian squares. By Proposition~\ref{prop:Gysin}(i), if $f_0$ and $f_1$ are lci and smoothable, the classes $\eta_{f_1}$ and $\eta_{f_2}\cdot\eta_{f_0}$ in $\pi_0\1^\BM(Z/X,\sL_{f_1})$ are equal, where $\eta_{f_2}\cdot\eta_{f_0}$ is the composite
\[
\Sigma^{\sL_{f_1}} \1_Z \simeq \Sigma^{f_2^*\sL_{f_0}} \Sigma^{\sL_{f_2}} \1_Z \xrightarrow{\eta_{f_2}} f_2^! \Sigma^{\sL_{f_0}} \1_Y \xrightarrow{\eta_{f_0}} f_2^!f_0^!\1_X \simeq f_1^!\1_X.
\]
These equalities form a section of the functor
\[
\Fun_{\mathrm{slci},\mathrm{slci},\mathrm{slci}}^{\mathrm{cart},\pitchfork}(\Delta^2,\Sch)^\op \to \tau_{\leq -1}\Spc,\quad (f_0,f_1,f_2)\mapsto \tau_{\leq -1}\Maps_{\1^\BM(Z/X,\sL_{f_1})}(\eta_{f_1},\eta_{f_2}\cdot\eta_{f_0}),
\]
This can be refined to a section of the functor
\[
\Fun_{M_0,M_1,M_2}^{\mathrm{cart},\pitchfork}(\Delta^2,\Sch)^\op \to \Spc,\quad (f_0,f_1,f_2)\mapsto \Maps_{\1^\BM(Z/X,\sL_{f_1})}(\eta_{f_1},\eta_{f_2}\cdot\eta_{f_0}),
\]
at least if each $M_i$ is either the class of regular immersions or that of smooth morphisms. One can reduce as in~\sssecref{sssec:eta-functoriality} to the case of regular closed immersions, where an explicit functorial homotopy $\eta_{f_1}\simeq \eta_{f_2}\cdot\eta_{f_0}$ is given by a double deformation to the normal cone \cite[3.2.19]{DJK}.


\section{Comparison of transfers}
\label{sec:transfers}

In this section we show that the framed transfers in cohomology provided by the motivic recognition principle are given by Gysin maps.
In \ssecref{ssec:fundtr}, we define the \emph{fundamental transfer} associated with a tangentially framed correspondence using Gysin maps. We then introduce in \ssecref{ssec:voevtr} the \emph{Voevodsky transfer} associated with an equationally framed correspondence, and we show that the Voevodsky transfer computes the fundamental transfer. Finally, in \ssecref{ssec:recogtr}, we show that the transfers obtained from the recognition principle agree with the Voevodsky transfer.

Throughout this section, we fix a base scheme $S$ and a motivic spectrum $E\in \SH(S)$. As explained in Remark~\ref{rem:E}, the spectrum $E$ is only used for readability purposes. 


\ssec{The fundamental transfer} 
\label{ssec:fundtr}

\sssec{}
Recall that a \emph{tangentially framed correspondence} between $S$-schemes $X$ and $Y$ is the data of a span
\begin{equation*}
  \begin{tikzcd}
     & Z \ar[swap]{ld}{f}\ar{rd}{h} & \\
    X &   & Y
  \end{tikzcd}
\end{equation*}
over $S$, where $f$ is finite syntomic, together with an equivalence $\tau\colon 0\simeq \sL_f$ in the $\infty$-groupoid $K(Z)$. We denote by $\Corr^\fr_S(X,Y)$ the $\infty$-groupoid of tangentially framed correspondences from $X$ to $Y$, defined as
\[
\Corr^\fr_S(X,Y) = \colim_{X\xleftarrow fZ \to Y} \Maps_{K(Z)}(0,\sL_f),
\]
where the colimit is taken over the groupoid of spans with $f$ finite syntomic.

\sssec{}\label{sssec:fund-tr}
Note that a finite syntomic morphism $f\colon Z\to X$ admits a \emph{canonical} factorization
\[
\begin{tikzcd}
	Z \ar{dr}[swap]{f} \ar[hook]{r} & \bV(f_*\sO_Z) \ar{d} \\
	& X,
\end{tikzcd}
\]
which we use to define the fundamental class $\eta_f\in \1^\BM(Z/X,\sL_f)$.

\begin{defn}\label{defn:fundamental-transfer}
	Let $X,Y\in\Sch_S$ and let $\alpha=(Z,f,h,\tau)$ be a tangentially framed correspondence from $X$ to $Y$ over $S$. For $E\in\SH(S)$, the \emph{fundamental transfer} $\tr_\eta(\alpha)\colon E(Y)\to E(X)$ is the composition
	\[
	E(Y) \xrightarrow{h^*} E(Z) \stackrel{\tau}{\simeq} E(Z,\sL_f) \xrightarrow{f_!} E(X).
	\]
\end{defn}

Using the functoriality of $f\mapsto \eta_f$ described in \sssecref{sssec:eta-functoriality}, we obtain a map
\[
\tr_\eta\colon \Corr^\fr_S(X,Y) \to \Maps(E(Y),E(X)), \quad \alpha\mapsto \tr_\eta(\alpha),
\]
natural in $(X,Y,E)\in \Sch_S^\op\times \Sch_S\times\SH(S)$.
If $X$ and $Y$ are smooth over $S$, then by the Yoneda lemma we obtain a map
\[
\Corr^\fr_S(X,Y) \to \Maps(\Sigma^\infty_\T X_+,\Sigma^\infty_\T Y_+),
\]
which we sometimes also denote by $\tr_\eta$.

\begin{rem}\label{rem:homological-gysin}
	Dually, a framed correspondence $\alpha$ as above also induces a map
	\begin{equation}\label{eqn:fund-transfer-M}
	M_S(X) \xrightarrow{f^!} M_S(Z,-\sL_f) \stackrel\tau\simeq M_S(Z) \xrightarrow{h_*} M_S(Y)
	\end{equation}
	in $\SH(S)$ (see~\sssecref{sssec:homological-motives} for the notation $M_S(X,\xi)$).
	By unpacking the definitions, it is easy to show that the natural transformation $E(-) \to \Maps(M_S(-),E)$ on $S$-schemes is also natural with respect to Gysin maps. In particular, if $X$ and $Y$ (but not necessarily $Z$) are cdh-locally smooth over $S$, applying $\Maps(-,E)$ to~\eqref{eqn:fund-transfer-M} yields the fundamental transfer $\tr_\eta(\alpha)\colon E(Y)\to E(X)$.
\end{rem}

\sssec{Example: the action of K-theory}
The $\infty$-groupoid $\Corr^\fr_S(S,S)$ contains $\Omega K(S)$ as a full subgroupoid. By construction, the composite
\[
\Omega K(S) \subset \Corr^\fr_S(S,S) \xrightarrow{\tr_\eta} \End(\1_S)
\]
is the action of $\Omega K(S)$ on $\1_S$ induced by the motivic J-homomorphism $K(S) \to \SH(S)$.

If $S$ is a regular semilocal scheme over a field of characteristic not $2$, $\pi_0\1_S(S)$ is isomorphic to the Grothendieck–Witt group $\GW(S)$ of nondegenerate symmetric bilinear forms over $S$ \cite[Lemma 10.12]{norms}. This isomorphism is such that the J-homomorphism
\[
\sO(S)^\times \simeq K_1(S) \to \pi_0\End(\1_S) \simeq \GW(S)
\]
sends a unit $a$ to the class $\langle a\rangle$ of the bilinear form $(x,y)\mapsto axy$.

\sssec{Example: finite étale transfers}
There is a canonical map $\Corr^\fet_S(X,Y) \to \Corr^\fr_S(X,Y)$, sending a span
\[
  \begin{tikzcd}
     & Z \ar[swap]{ld}{f}\ar{rd} & \\
    X &   & Y
  \end{tikzcd}
\]
with $f$ finite étale to the same span equipped with the canonical trivialization of $\sL_f$. 

If $S$ is a regular semilocal scheme over a field of characteristic not $2$ and $\alpha$ is the finite étale correspondence $S\xleftarrow f T\xrightarrow{\id} T$, the transfer
\[
\tr_\eta(\alpha)\colon \GW(T) \to \GW(S)
\]
is the Scharlau transfer associated with the canonical trace $\mathrm{Tr}_{T/S}\colon \sO(T)\to \sO(S)$. Indeed, one is reduced to the case of a field extension using the Gersten resolution for Grothendieck–Witt groups \cite[Theorem 100]{BalmerWitt}, and in that case the claim was proved in \cite[\sectsign 5]{HoyoisGLV}.

\sssec{The oriented case}
Suppose given a retraction diagram
\[
E \xrightarrow\iota \MGL\otimes E \xrightarrow\rho E, \quad \rho\circ\iota \simeq \id_E,
\]
in $\SH(S)$; such a diagram exists if $E$ is an $\MGL$-module in the homotopy category $\ho\SH(S)$, and it is given if $E$ is an $\MGL$-module in $\SH(S)$. Then the fundamental transfers in $E$-cohomology are independent of the tangential framings.
More precisely, given $X,Y\in\Sch_S$, there is a canonical factorization
\[
\begin{tikzcd}
	\Corr^\fr_S(X,Y) \ar{d}[swap]{\text{forget}} \ar{r}{\tr_\eta} & \Maps(E(Y),E(X)). \\
	\Corr^\fsyn_S(X,Y) \ar[dashed]{ur} & 
\end{tikzcd}
\] 
This follows at once from the fact that the $\MGL$-linearized J-homomorphism
\[
K(S) \to \SH(S), \quad \xi \mapsto \Sigma^\xi\MGL
\]
factors through the rank map $\rk\colon K\to \Z$ \cite[\sectsign 16.2]{norms}.



\ssec{The Voevodsky transfer} 
\label{ssec:voevtr}

\sssec{}

\newcommand{\can}{\mathrm{can}}

Let $X$ and $Y$ be $S$-schemes and let $\alpha\in\Corr^{\efr,n}_S(X,Y)$ be an equationally framed correspondence of level $n$ from $X$ to $Y$ \cite[Definition 2.1.2]{EHKSY1}. We display $\alpha$ as the diagram
\begin{equation}\label{eq:efr-diagram}
\begin{tikzcd}
  & U \ar{dr}{(\phi, g)} \ar{dl}[swap]{u} &  \\
  \A^n_X \ar{d}[swap]{\pi} & Z \ar[hook]{u}{i} \ar{dl}{f} \ar[hook']{l}[swap]{i} \ar{dr}[swap]{h} & \A^n \times Y \\ 
 X & & Y, \ar[hook]{u}[swap]{0}
\end{tikzcd}
\end{equation}
where $f$ is finite, $u$ is an étale neighborhood of $Z$ in $\A^n_X$, $0$ is the zero section, and the right-hand square is Cartesian. 

We will denote by $\P^{\times n}$ the $n$-fold product $(\P^1)^{\times n}$, regarded as a compactification of $\A^n$, and by $\partial\P^{\times n}\subset \P^{\times n}$ the complementary reduced closed subscheme that is the union of the ``faces'' $\P^{\times i-1}\times\{\infty\}\times\P^{\times n-i}$:
\begin{equation*}
\begin{tikzcd}
	\A^n_X \ar[hook]{r}{j} \ar{dr}[swap]{\pi} & \P_X^{\times n} \ar{d}{\hat\pi} & \partial\P_X^{\times n}. \ar[hook']{l}[swap]{\infty} \ar{dl} \\
	& X & 
\end{tikzcd}
\end{equation*}

\begin{defn}\label{defn:Voevodsky-transfer}
	Let $X,Y\in\Sch_S$ and let $\alpha\in \Corr^{\efr,n}_S(X,Y)$ be the equationally framed correspondence~\eqref{eq:efr-diagram}. For $E\in\SH(S)$, the \emph{Voevodsky transfer} $\tr_V(\alpha)\colon E(Y)\to E(X)$ is the composition
	\begin{multline*}
	E(Y) \simeq E_Y(\A^n_Y,\sO^n) \xrightarrow{(\phi,g)^*} E_Z(U,\sO^n) \simeq E_Z(\P_X^{\times n},\sO^n)
	\\
	 \to E_{\A^n_X} (\P_X^{\times n},\sO^n) \xleftarrow{\simeq} E_X(\P_X^{\times n},\sO^n) \simeq E(X).
	\end{multline*}
\end{defn}

Here, the maps $E_Z(\P_X^{\times n},\sO^n)\to E_{\A^n_X}(\P_X^{\times n},\sO^n) \leftarrow E_X(\P_X^{\times n},\sO^n)$ are instances of~\eqref{eqn:forget-support}.
To see that the latter is an equivalence, first note that it fits in the diagram
\[
\begin{tikzcd}
	E_X(\P_X^{\times n},\sO^n) \ar{r} \ar[dashed]{d} & E(\P_X^{\times n},\sO^n) \ar{r} \ar{d}{\id} & E(\P_X^{\times n}\setminus 0_X,\sO^n) \ar{d} \\
	E_{\A^n_X}(\P_X^{\times n},\sO^n) \ar{r} &  E(\P_X^{\times n},\sO^n)\ar{r}& E(\partial\P_X^{\times n},\sO^n),
\end{tikzcd}
\]
where the rows are the fiber sequences \eqref{eqn:E_Z(X)} and~\eqref{eqn:E(X/Z)}. By~\ref{sssec: descent properties}, the claim then follows from the following lemma:

\begin{lem}\label{lem:boundary}
	The inclusion $\partial\P^{\times n}\subset \P^{\times n}\setminus 0$ is an $\A^1$-cdh-equivalence over $\Spec\Z$.
\end{lem}

\begin{proof}
We consider the commutative square of inclusions in $\Pre(\Sch)$:
\[
\begin{tikzcd}
	\bigcup_{i=1}^n \left( \P^{\times i-1} \times \{\infty\}\times \P^{\times n-i}\right) \ar{r} \ar{d} & \partial\P^{\times n} \ar{d} \\
	\bigcup_{i=1}^n \left( \P^{\times i-1} \times (\P^1\setminus 0) \times \P^{\times n-i}\right) \ar{r} & \P^{\times n}\setminus 0.
\end{tikzcd}
\]
The upper horizontal map is a covering sieve in the closed topology, the lower horizontal map is a covering sieve in the open topology, and the left vertical map is the colimit of an $n$-dimensional cube of $\A^1$-homotopy equivalences. In particular, these three maps are $\A^1$-cdh-equivalences, hence so is the right vertical map.
\end{proof}

\begin{rem}
	In general, $U$ is an algebraic space and not a scheme, but this does not matter. Indeed, the inclusion of schemes into (Zariski-locally quasi-separated) algebraic spaces induces an equivalence between the $\infty$-categories of Nisnevich sheaves, by \cite[Proposition 5.7.6]{GrusonRaynaud}. As a result, we may tacitly extend any Nisnevich sheaf, such as $E(-)$ or $\SH(-)$, to algebraic spaces. However, we can assume that $U$ is a scheme in many cases \cite[Lemma A.1.2(iv)]{EHKSY1}.
\end{rem}

\sssec{}
By Voevodsky's Lemma \cite[Corollary A.1.7]{EHKSY1}, the equationally framed correspondence $\alpha$ is equivalently a morphism of pointed presheaves $\Sigma^n_\P X_+ \to L_\nis \Sigma^n_\T Y_+$. Explicitly, it is given by the following zig-zag in $\Pre(\Sch_S)_*$:
\begin{equation}\label{eqn:voev-lemma}
	 \Sigma^{n}_{\P} X_+ = \frac{\P_X^{\times n}}{\tilde\partial\P_X^{\times n}} \rightarrow \frac{\P_X^{\times n}}{\P_X^{\times n} \setminus Z}  \xleftarrow{u}  \frac{U}{U \setminus Z} \xrightarrow{(\phi,g)} \frac{\A^n_Y}{\A^n_Y \setminus 0} \leftarrow \Sigma^n_\T Y_+.
\end{equation}
Here, $\tilde\partial\P^{\times n}\subset \P^{\times n}$ is the subpresheaf defined as the union
\[
\tilde\partial\P^{\times n}=\bigcup_{i=1}^n \left(\P^{\times i-1}\times \{\infty\}\times \P^{\times n-i}\right),
\] 
the first map is the collapse map, and the wrong-way maps are Nisnevich-local equivalences. The Voevodsky transfer $\tr_V(\alpha)\colon E(Y)\to E(X)$ is then equivalent to applying (the right Kan extension of) the functor $E(-,\sO^n)$ to the composite~\eqref{eqn:voev-lemma}. 

In particular, if $X$ and $Y$ are smooth over $S$, then~\eqref{eqn:voev-lemma} induces a morphism $\Sigma^\infty_\T X_+ \to \Sigma^\infty_\T Y_+$ in $\SH(S)$, which gives the Voevodsky transfer upon applying $\Maps(-,E)$.

\sssec{} For every $n\geq 0$, Definition~\ref{defn:Voevodsky-transfer} gives a map
\[
\tr_V\colon \Corr^{\efr,n}_S(X,Y) \to \Maps(E(Y),E(X)).
\]
Using for instance~\eqref{eqn:voev-lemma}, it is clear that this map is natural in $(X,Y)\in\Sch_S^\op\times\Sch_S$. Moreover, by \cite[Remark 2.1.6]{EHKSY1}, the triangles
\[
\begin{tikzcd}
	\Corr^{\efr,n}_S(X,Y) \ar{r}{\tr_V} \ar{d}[swap]{\sigma} & \Maps(E(Y),E(X)) \\
	\Corr^{\efr,n+1}_S(X,Y) \ar{ur}[swap]{\tr_V} & 
\end{tikzcd}
\]
naturally commute (here $\sigma$ is the suspension morphism~\cite[2.1.4]{EHKSY1}). Passing to the colimit gives a natural map
\[
\tr_V\colon \Corr^{\efr}_S(X,Y) \to \Maps(E(Y),E(X)).
\]

If we let $X$ vary, note that $\Maps(E(Y),E(-))\colon \Sch_S^\op \to \Spc$ is an $\A^1$-invariant cdh sheaf.
If $Y$ is \emph{smooth} over $S$, the forgetful map $\Corr^{\efr}_S(-,Y) \to \Corr^{\fr}_S(-,Y)$ is a motivic equivalence \cite[Corollary 2.3.27]{EHKSY1}. In that case, therefore, the Voevodsky transfer factors through the $\infty$-groupoid of tangentially framed correspondences, inducing a morphism
\begin{equation}\label{eq:voev-transfer-fr}
\tr_V\colon \Corr^{\fr}_S(X,Y) \to \Maps(E(Y),E(X)).
\end{equation}

\sssec{}
We now prove that the Voevodsky transfer agrees with the fundamental transfer of Definition~\ref{defn:fundamental-transfer}.

\begin{thm} \label{thm:main-comparison}
Let $X,Y\in\Sch_S$ and $E\in\SH(S)$. Then the triangle
\[
\begin{tikzcd}
	\Corr^\efr_S(X,Y) \ar{d}[swap]{\textnormal{forget}} \ar{r}{\tr_V} & \Maps(E(Y),E(X)) \\
	\Corr^\fr_S(X,Y) \ar{ur}[swap]{\tr_\eta} & 
\end{tikzcd}
\]
commutes, naturally in $E$, $X$, and $Y$.
In particular, if $Y$ is smooth over $S$, then the Voevodsky transfer~\eqref{eq:voev-transfer-fr} coincides with the fundamental transfer.
\end{thm}

\begin{proof}
	Let $\alpha$ be an equationally framed correspondence as in~\eqref{eq:efr-diagram}, and let $\tau\colon 0\simeq\sL_f$ be the induced trivialization in $K(Z)$. We must show that the following diagram commutes:
	\begin{equation*}
	\begin{tikzcd}
		E(Y) \ar{r}{\simeq} \ar{d}[swap]{h^*} & E_Y(\A^n_Y,\sO^n)  \ar{d}{(\phi,g)^*} \\
		E(Z) \ar{d}{\simeq}[swap]{\tau} & E_Z(U,\sO^n) \\
		E(Z,\sL_f)  \ar{dd}[swap]{f_!} & E_Z(\P^{\times n}_X,\sO^n)  \ar{u}[swap]{\simeq}  \ar{d}  \\
		 & E_{\A^n_X}(\P^{\times n}_X,\sO^n) \\
		E(X)  \ar{r}{\simeq} & E_X(\P^{\times n}_X,\sO^n).  \ar{u}[swap]{\simeq}
	\end{tikzcd}
	\end{equation*}
	To do so, we subdivide this diagram as follows:
	\begin{equation}\label{eq:transfer-comparison}
		\begin{tikzcd}
			E(Y) \ar{rr}{0_!}[swap]{\simeq} \arrow[drr, phantom, "\scriptstyle\text{(1)}"] \ar{d}[swap]{h^*} & & E_Y(\A^n_Y,\sO^n)  \ar{d}{(\phi,g)^*} \\
			E(Z) \ar{d}{\simeq}[swap]{\tau} \ar{r}{i_!} & E_Z(U,-\sL_i) \ar{r}{\simeq} \arrow[dr, phantom, near end, "\scriptstyle\text{(3)}"] &  E_Z(U,\sO^n) \ar[bend right]{d}[swap]{u_!} \\
			E(Z,\sL_f) \ar{dr}{i_!} \ar{dd}[swap]{f_!} \ar{rru}{i_!} \ar{rr}{i_!} & \arrow[ul, phantom, near end, "\scriptstyle(2)"] & E_Z(\P^{\times n}_X,\sO^n)  \ar{u}[swap]{\simeq}  \ar{d} \ar{dl}  \\
			\arrow[r, phantom, "\scriptstyle\text{(4)}"] & E(\P_X^{\times n},\sL_{\hat\pi}) \ar{dl}{\hat\pi_!} & E_{\A^n_X}(\P^{\times n}_X,\sO^n) \ar{l} \\
			E(X)  \ar{rr}{0_!}[swap]{\simeq} & \arrow[u, phantom, "\scriptstyle\text{(5)}"] & E_X(\P^{\times n}_X,\sO^n).  \ar{u}[swap]{\simeq} \ar{ul}
		\end{tikzcd}
	\end{equation}
The rectangle (1) commutes by the base change property of Gysin maps (Proposition~\ref{prop:Gysin}(ii)) applied to the Cartesian square
\[
\begin{tikzcd}
	Z \ar{r}{h} \ar{d}[swap]{i} & Y \ar{d}{0} \\
	U \ar{r}{(\phi,g)} & \A^n\times Y,
\end{tikzcd}
\]
which is tor-independent since $i$ is a regular immersion of codimension $n$.
Thus, the unnamed equivalence in (1) is induced by the isomorphism $\sN_i\simeq h^*(\sN_0)\simeq\sO^n$. This isomorphism also induces the trivialization $\tau$, whence the commutativiy of the square (2). The triangles (3), (4), and (5) all commute by the compatibility of Gysin maps with composition (Proposition~\ref{prop:Gysin}(i)), where the commutativity of (5) means that going around starting from the lower left corner gives the identity.

To conclude the proof, we must show that the diagram~\eqref{eq:transfer-comparison} can be promoted to a functor of the triple $(X,Y,\alpha)$. This follows from the functoriality properties of Gysin maps discussed in \sssecref{sssec:eta-functoriality} and \sssecref{sssec:eta-composition}. For the triangle (4), we must recall that the fundamental class $\eta_f$ was defined using the canonical factorization $Z\hook \bV(f_*\sO_Z)\to X$. The commutativity of (4) can be made functorial using the triangles
\[
\begin{tikzcd}
	& Z \ar[hook']{dl} \ar[hook]{dr} \ar[hook]{d} & \\
	\bV(f_*\sO_Z) \ar{dr} & \bV(f_*\sO_Z) \times \P^{\times n} \ar{l} \ar{r} \ar{d} & \P_X^{\times n} \ar{dl} \\
	& X &
\end{tikzcd}
\]
in which the upper three maps are regular closed immersions and the other five are smooth.
This concludes the proof of the theorem.
\end{proof}


\ssec{The transfer from the recognition principle}
\label{ssec:recogtr}

\sssec{}
Recall that there is an $\infty$-category $\Span^\fr(\Sm_S)$ whose objects are smooth $S$-schemes and whose mapping spaces are the $\infty$-groupoids $\Corr^\fr_S(X,Y)$, which gives rise to the $\infty$-category $\SH^\fr(S)$ of \emph{framed motivic spectra} \cite[\sectsign 3]{EHKSY1}. The ``graph'' functor
\[
\gamma\colon \Sm_{S+} \to \Span^\fr(\Sm_S), \quad (f\colon X_+\to Y_+) \mapsto (X\hookleftarrow f^{-1}(Y) \xrightarrow{f} Y),
\]
induces an adjunction
\[
\gamma^*: \SH(S) \rightleftarrows \SH^\fr(S): \gamma_*
\]
such that the following square commutes:
\[
\begin{tikzcd}
	\Sm_{S+} \ar{r}{\gamma} \ar{d}[swap]{\Sigma^\infty_\T} & \Span^\fr(\Sm_S) \ar{d}{\Sigma^\infty_{\T,\fr}} \\
	\SH(S) \ar{r}{\gamma^*} & \SH^\fr(S).
\end{tikzcd}
\]

By the reconstruction theorem \cite[Theorem 16]{Hoyois:2018aa}, the functor $\gamma^*\colon \SH(S)\to \SH^\fr(S)$ is an equivalence of $\infty$-categories. It follows that $E$-cohomology of smooth $S$-schemes acquires canonical framed transfers:
\[
\begin{tikzcd}
	\Sm_S^\op \ar{r}{E(-)} \ar{d}[swap]{\gamma} & \Spc. \\
	\Span^\fr(\Sm_S)^\op \ar[dashed]{ur}
\end{tikzcd}
\]
The goal of this section is to show that these transfers coincide with the Voevodsky transfers, hence with the fundamental transfers.

\sssec{}

We need a technical preliminary result, which we formulate in a more general context.

Let $C$ be a presentably symmetric monoidal $\infty$-category and let $\T\in C$ be an object. For any presentable $C$-module $M$, we have an adjunction
\[
\Sigma_\T: M\rightleftarrows M:\Omega_\T
\]
where $\Sigma_\T=\T\otimes(-)$. We can then form the diagram
\[
\N\times\N \to \Fun(M,M), \quad (m,n)\mapsto \Omega^n_\T\Omega^m_\T\Sigma^m_\T\Sigma^n_\T,
\]
where the transition maps use the unit transformation $\id\to\Omega_\T\Sigma_\T$ and (in one direction) the cyclic permutations $\Sigma_\T\Sigma^m_\T\simeq \Sigma^m_\T\Sigma_\T$ and $\Omega^m_\T \Omega_\T\simeq \Omega_\T \Omega^m_\T$.

We denote by $\Spt_\T(M)$ the $\infty$-category of $\T$-spectra in $M$, defined as the limit
\[
\Spt_\T(M) = \lim\left(\dotsb \to M \xrightarrow{\Omega_\T} M \xrightarrow{\Omega_\T} M\right).
\]
We have an adjunction
\[
\Sigma^\infty_\T: M \rightleftarrows \Spt_\T(M): \Omega^\infty_\T,
\]
where $\Omega^\infty_\T$ is the projection to the last copy of $M$. If $\Omega_\T$ preserves sequential colimits, then
\[
\Omega_\T^\infty\Sigma^\infty_\T \simeq \colim_n \Omega_\T^n\Sigma_\T^n.
\]

\begin{lem}\label{lem:T-spectra-business}
	With the above notation, suppose that $\Omega_\T\colon M\to M$ preserves sequential colimits and that the cyclic permutation of $\T^{\otimes n}$ is homotopic to the identity for some $n\geq 2$.
	Then the natural transformations
	\begin{gather*}
		\Omega^\infty_\T\Sigma^\infty_\T \simeq \colim_p \Omega^p_\T\Sigma^p_\T \to \colim_{p,q} \Omega^p_\T\Omega^q_\T\Sigma^q_\T\Sigma^p_\T,\\
		\Omega^\infty_\T\Sigma^\infty_\T \simeq \colim_q \Omega^q_\T\Sigma^q_\T \to \colim_{p,q} \Omega^p_\T\Omega^q_\T\Sigma^q_\T\Sigma^p_\T
	\end{gather*}
	between endofunctors of $M$ are homotopic equivalences.
\end{lem}

\begin{proof}
	We have a commutative diagram
	\[
	\begin{tikzcd}
		M \ar{r}{\Sigma^\infty_\T} & \Spt_\T(M) \arrow[r, shift left, "F_0"]
	\arrow[r, shift right, swap, "F_1"] & \Spt_\T(\Spt_\T(M)),
	\end{tikzcd}
	\]
	where $F_0=\Sigma^\infty_\T$ and $F_1=\Spt_\T(\Sigma^\infty_\T)$.
	Let $G_i$ be the right adjoint to $F_i$ and $u_i\colon \id\to G_iF_i$ the unit transformation.
	Then the given natural transformations are $\Omega^\infty_\T u_i \Sigma^\infty_\T$ for $i=0,1$.
	
	By the assumption on $\T$, the functor $\Spt_\T(-)$ is a left localization of the $\infty$-category of presentable $C$-modules \cite[Corollary 2.22]{Robalo}. This implies that $F_0$ and $F_1$ are equivalences and moreover that there is a natural equivalence $\alpha\colon F_0\simeq F_1$ such that $\alpha \Sigma^\infty_\T$ is the identity. In particular, the unit transformations $u_i$ are equivalences, and $\alpha$ and its mate $G_0\simeq G_1$ give the desired homotopy.
\end{proof}

\sssec{} We now prove that the Voevodsky transfer coincides with the transfer coming from the reconstruction theorem.
For $X,Y\in\Sm_S$, we can regard the Voevodsky transfer as a map
\begin{equation}\label{eqn:efr-to-SH}
\tr_V\colon \Corr^\efr_S(X,Y) \to \Maps_{\SH(S)}(\Sigma^\infty_\T X_+, \Sigma^\infty_\T Y_+),
\end{equation}
sending the equationally framed correspondence~\eqref{eq:efr-diagram} to the composite~\eqref{eqn:voev-lemma}. Since $Y$ is smooth, this map factors through the $\infty$-groupoid $\Corr^\fr_S(X,Y)$.

\begin{thm}\label{thm:comparison2}
	For $X,Y\in\Sm_S$, the following diagram commutes, naturally in $X$ and $Y$:
	\[
	\begin{tikzcd}
		\Corr^\fr_S(X,Y) \ar{r}{\tr_V} \ar{dr}[swap]{\Sigma^\infty_{\T,\fr}} & \Maps(\Sigma^\infty_\T X_+,\Sigma^\infty_\T Y_+) \ar{d}{\gamma^*}[swap]{\simeq} \\
		& \Maps(\gamma^*\Sigma^\infty_\T X_+,\gamma^*\Sigma^\infty_\T Y_+).
	\end{tikzcd}
	\]
\end{thm}

\begin{proof}
	Recall the endofunctors $\h^{\efr,n}_S$ and $\h^\efr_S$ of $\Pre_\Sigma(\Sm_S)_*$ defined in \cite[2.1.10]{EHKSY1}.
	By Voevodsky's Lemma, there is a canonical transformation
		\begin{equation*}\label{eqn:efr-to-OmegaSigma}
		 \h^{\efr,n}_S\to \Omega^n_\P L_\nis \Sigma^n_\T\colon \Pre_\Sigma(\Sm_S)_* \to \Pre_\Sigma(\Sm_S)_*,
		\end{equation*}
		extending an equivalence on representables. Composing with $L_\nis\to L_\mot$ and taking the colimit over $n$, we obtain a transformation
		\[
		\h^\efr_S \to \Omega^\infty_\T\Sigma^\infty_\T L_\mot,
		\]
		which extends \eqref{eqn:efr-to-SH} to pointed presheaves.
	We will prove more generally that the following diagram of endofunctors of $\Pre_\Sigma(\Sm_S)_*$ commutes:
	\begin{equation}\label{eqn:comparison}
		\begin{tikzcd}
			 \h^\efr_S \ar{r} \ar[swap]{d} & \Omega^\infty_\T\Sigma^\infty_\T L_\mot \ar{d}{\text{unit}} \\
			 \h^\fr_S \ar[swap]{d} & \Omega^\infty_\T \gamma_*\gamma^*\Sigma^\infty_\T L_\mot \ar{d}{\simeq} \\
			\gamma_*\gamma^* L_\mot \ar{r}{\text{unit}} &  \gamma_* \Omega^\infty_{\T,\fr}\Sigma^\infty_{\T,\fr}\gamma^* L_\mot.
		\end{tikzcd}
	\end{equation}
	Since $\h^\efr_S$ is by definition the left Kan extension of its restriction to $\Sm_{S+}$, it suffices to show that the diagram commutes on $\Sm_{S+}$.  Recall that the unit map $\id\to \gamma_*\gamma^* \simeq \h^\fr_S$ factors as
	\[
	\id \to \h^\efr_S \to \h^\fr_S,
	\]
	where the second map is a motivic equivalence \cite[Corollary 2.3.27]{EHKSY1}.

Consider the two commuting squares
	\[
	\begin{tikzcd}[column sep=2em]
		\colim_{n} \h^{\efr,n}_S \ar[shift left]{r}{q=0} \ar[swap,shift right]{r}{p=0} \ar{d} & 
		\colim_{p,q}\Omega^p_\P L_\mot \h^{\efr,q}_S\Sigma^p_\T \ar{d} \\
		\colim_n \Omega_\P^n L_\mot \Sigma_\T^n \ar[shift left]{r}{q=0} \ar[swap,shift right]{r}{p=0} &
		\colim_{p,q}\Omega^p_\P \Omega^q_\P L_\mot \Sigma^q_\T\Sigma^p_\T,
	\end{tikzcd}
	\]
	where the vertical maps are induced by~\eqref{eqn:efr-to-SH}.
	The two upper horizontal maps can be identified with the two composites in the diagram \eqref{eqn:comparison}.
	On the other hand, the two lower horizontal maps are equivalent by Lemma~\ref{lem:T-spectra-business}.
	It therefore suffices to show that the right vertical map is an equivalence.
	
	We have a commuting triangle
	\[
	\begin{tikzcd}
		\colim_p \Omega^p_\P L_\mot \Sigma^p_\T \ar{d} \ar{dr} & \\
		\colim_{p,q}\Omega^p_\P L_\mot \h^{\efr,q}_S\Sigma^p_\T \ar{r} & \colim_{p,q}\Omega^p_\P \Omega^q_\P L_\mot \Sigma^q_\T\Sigma^p_\T.
	\end{tikzcd}
	\]
	The diagonal map is an equivalence by Lemma~\ref{lem:T-spectra-business}. The vertical map is the unit map $\Omega^\infty_\T\Sigma^\infty_\T \to \Omega^\infty_\T\gamma_*\gamma^*\Sigma^\infty_\T$, which is an equivalence since $\gamma^*\colon \SH(S)\to \SH^\fr(S)$ is fully faithful. Hence, the bottom horizontal map is an equivalence, as desired.
\end{proof}

\begin{cor} \label{cor:efr-vs-maps}
	Let $k$ be a perfect field and let $\sF\in\Pre_\Sigma(\Sm_k)_*$. Then the map \eqref{eqn:efr-to-SH} induces an equivalence
	\[
	L_\zar(\Lhtp \h^\efr_k(\sF))^\gp \simeq \Omega^\infty_\T\Sigma^\infty_\T \sF\colon \Sm_k^\op \to \CMon^\gp(\Spc).
	\]
\end{cor}

Note that we already have such an equivalence by \cite[Corollary 3.5.16]{EHKSY1}. The point of this corollary is that this equivalence is induced by Voevodsky's Lemma, as one would expect.

\begin{proof}
	If we plug $\sF$ in the diagram~\eqref{eqn:comparison}, we obtain a commutative square of $\Fin_*$-objects, all of which are $\Einfty$-objects except the top left corner. 
	The upper left vertical map is a motivic equivalence by \cite[Corollary 2.3.27]{EHKSY1}, and so is the lower left vertical map since $\gamma_*\gamma^*$ preserves motivic equivalences \cite[Proposition 3.2.14]{EHKSY1}. The right vertical map is an equivalence by \cite[Theorem 3.5.12]{EHKSY1}, and the lower horizontal map is group completion by \cite[Corollary 3.5.10]{EHKSY1}. The commutativity of the diagram shows that \eqref{eqn:efr-to-SH} induces an equivalence
	\[
	(L_\mot\h^\efr_k(\sF))^\gp \simeq \Maps_{\SH(k)}(\Sigma^\infty_\T (-)_+,\Sigma^\infty_\T \sF).
	\]
	On the other hand, the canonical map
	\[
	L_\zar(\Lhtp \h^\efr_k(\sF))^\gp \to (L_\mot\h^\efr_k(\sF))^\gp
	\]
	is an equivalence since the left-hand side is already Nisnevich-local and $\A^1$-invariant \cite[Corollary 3.5.16]{EHKSY1}.
\end{proof}

\sssec{} Combining Theorems \ref{thm:main-comparison} and \ref{thm:comparison2} we obtain:

\begin{thm}\label{thm:fund-vs-recog}
	For $X,Y\in\Sm_S$, the following diagram commutes, naturally in $X$ and $Y$:
	\[
	\begin{tikzcd}
		\Corr^\fr_S(X,Y) \ar{r}{\tr_\eta} \ar{dr}[swap]{\Sigma^\infty_{\T,\fr}} & \Maps(\Sigma^\infty_\T X_+,\Sigma^\infty_\T Y_+) \ar{d}{\gamma^*}[swap]{\simeq} \\
		& \Maps(\gamma^*\Sigma^\infty_\T X_+,\gamma^*\Sigma^\infty_\T Y_+).
	\end{tikzcd}
	\]
\end{thm}

\begin{cor}\label{cor:action-of-K}
	Let $X$ be a smooth $S$-scheme. Then the
	motivic J-homomorphism $\Omega K(X) \to \End(\Sigma^\infty_\T X_+)$ coincides with the composition
	\[
	\Omega K(X) \subset \Corr_S^\fr(X,X) \xrightarrow{\Sigma^\infty_{\T,\fr}} \End(\Sigma^\infty_{\T,\fr} X) \simeq \End(\Sigma^\infty_\T X_+).
	\]
\end{cor}

\begin{proof}
	This follows immediately from Theorem~\ref{thm:fund-vs-recog}.
\end{proof}


\section{Finite correspondences for motivic ring spectra}
\label{sec:fin-corr}

In this section we introduce finite $R$-correspondences for a motivic ring spectrum $R$, generalizing the finite correspondences of Voevodsky and the finite Milnor–Witt correspondences of Calmès and Fasel.
In \ssecref{ssec:borel-moore}, we construct the homotopy category $\ho\Span^R(\Sm_S)$ of finite $R$-correspondences between smooth $S$-schemes, together with a functor to the homotopy category of $R$-modules. In \ssecref{ssec:functor}, we construct a functor from the category of (tangentially) framed correspondences to that of finite $R$-correspondences and compare it with the free $R$-module functor. Finally, in \ssecref{ssec:VMW}, we compare our constructions with those of Voevodsky and of Calmès–Fasel.

Throughout this section, $S$ is a fixed base scheme. All $S$-schemes are assumed to be separated.


\ssec{The category of finite $R$-correspondences} \label{ssec:borel-moore} 
Given an associative ring spectrum $R \in \SH(S)$, we will construct an $\ho\Spc$-enriched category $\ho\Span^R(\Sm_S)$ of \emph{finite $R$-correspondences} between smooth $S$-schemes.

To motivate our construction, recall that a morphism from $X$ to $Y$ in Voevodsky's category of finite correspondences over a regular scheme $S$ is an element of the free abelian group generated by integral closed subschemes $Z \subset X \times_S Y$ that are finite and surjective over a component of $X$. Alternatively, we can think of a morphism in this category as a reduced closed subscheme $Z \subset X \times_S Y$, each of whose irreducible components is finite and surjective over a component of $X$ and labeled by an integer. The category $\ho\Span^R(\Sm_S)$ will admit a similar description, but with integers replaced by Borel–Moore $R$-homology classes of $Z$ over $X$.

\sssec{} Let $S$ be a scheme and $R\in \SH(S)$.
For separated $S$-schemes $X,Y\in\Sch_S$, define
\[
\Corr^R_S(X,Y) =  \colim_{Z \subset X \times_S Y} R^{\BM}(Z/X),
\]
where the colimit is taken over the filtered poset of reduced subschemes $Z \subset X \times_S Y$ that are finite and universally open\footnote{equivalently, when $X$ is Noetherian, universally equidimensional \cite[Proposition 2.1.7]{sv-cycleschow}} over $X$. To form this colimit, we use the covariant functoriality of Borel–Moore homology with respect to proper morphisms (see \sssecref{sssec:functoriality}).

\sssec{}\label{sssec:GammaR}
Suppose now that $R\in\SH(S)$ is a homotopy associative ring spectrum, i.e., $R$ is an associative algebra in the homotopy category $\h\SH(S)$.
There is a map
\[
\Gamma^R\colon \Maps_S(X, Y) \to \Corr^R_S(X,Y)
\]
sending an $S$-morphism $f\colon X\to Y$ to its graph $\Gamma_f\subset X\times_SY$ labeled by the unit element $1\in R^\BM(\Gamma_f/X) \simeq \Maps(\1_X,R_X)$.

\sssec{}\label{sssec:CorrR-composition}
For $X,Y,T\in\Sch_S$, we define a composition law
\begin{equation}\label{eq:composition in Corr^R}
\circ \colon \Corr^R_S(X, Y) \times \Corr^R_S(Y, T) \to \Corr^R_S(X, T)
\end{equation}
as follows.
For any closed subschemes $Z \subset X \times_S Y$ and $Z' \subset Y \times_S T$, finite and universally open over $X$ and $Y$ respectively,
we consider the $2$-span
 \begin{equation*}
 \begin{tikzcd}[column sep={4.5em,between origins}]
 & & Z'' \ar{dr}{g'} \ar{dl}[swap]{h'} &  & \\
 & Z \ar{dr}{g} \ar{dl}[swap]{f} & &  Z' \ar{dl}[swap]{h} \ar{dr}{k}  & \\
 X & & Y & & T\rlap.
 \end{tikzcd}
 \end{equation*}
Let $p\colon Z''\to X\times_S T$ be the induced map and let $Z'\circ Z\subset X\times_ST$ be its reduced image.
It is clear that $Z'\circ Z$ is finite and universally open over $X$. We define the pairing
\[
\theta^\BM\colon R^{\BM}(Z/X) \times R^{\BM}(Z'/Y) \to R^{\BM}(Z'\circ Z/X)
\]
as the composition
\begin{multline*}
R^{\BM}(Z/X) \times R^{\BM}(Z'/Y) \xrightarrow{\id\times g^*} R^{\BM}(Z/X) \times R^{\BM}(Z''/Z)
         \\ \xrightarrow{\mu^\BM} R^{\BM}(Z''/X) \xrightarrow{p_*} R^{\BM}(Z'\circ Z/X).
\end{multline*}
 Here, $\mu^\BM$ is the composition product (which uses the ring structure on $R$, see \sssecref{sssec:products}(3)), and $p_*$ is the proper pushforward in Borel--Moore homology. More succinctly,
 \[
 \theta^\BM(x,y)=p_*(g^*(y)\cdot x).
 \] 
The map~\eqref{eq:composition in Corr^R} is then the filtered colimit over $Z$ and $Z'$ of the maps $\theta^\BM$. 

\begin{lem}\label{lem:associativity}
The composition law~\eqref{eq:composition in Corr^R} is unital and associative up to homotopy, with identity $\Gamma^R(\id_S)\in\Corr^R_S(X,X)$.
\end{lem}

\begin{proof}
Let $X_1$, $X_2$, $X_3$ and $X_4$ be smooth $S$-schemes, and let $Z_{i,i+1} \subset X_{i} \fibprod_S X_{i+1}$ be reduced subschemes, finite and universally open over $X_i$. 
Let $x_i \in R^{\BM}(Z_{i,i+1}/X_i)$ for each $1\le i\le 3$.
Consider the diagram
  \begin{equation*}
    \begin{tikzcd}[column sep={4.5em,between origins}]
      & & & Z_{1234} \ar{ld}\ar{rd} & & &
    \\
      & & Z_{123}\ar{rd}\ar{ld} & & Z_{234}\ar{rd}\ar{ld} & &
    \\
      & Z_{12}\ar{rd}\ar{ld} & & Z_{23}\ar{rd}\ar{ld} & & Z_{34}\ar{rd}\ar{ld} &
    \\
      X_1 & & X_2 & & X_3 & & X_4
    \end{tikzcd}
  \end{equation*}
where the squares are Cartesian.
It suffices to note that the two possible ways of composing the $x_i$'s are both equal to
  \begin{equation*}
    p_*(h^*(x_3) \cdot g^*(x_2) \cdot x_1),
  \end{equation*}
where $g\colon Z_{12}\to X_2$, $h\colon Z_{123}\to X_3$, and $p\colon Z_{1234} \to Z_{34}\circ Z_{23}\circ Z_{12}$.
This follows directly from the properties of the composition product listed in \cite[1.2.8]{deglise2017bivariant}. The fact that $\Gamma^R(\id_X)$ is the identity is trivial.
\end{proof}

\sssec{The category of finite $R$-correspondences}
In view of Lemma~\ref{lem:associativity}, we can define a category $\ho\Span^R(\Sch_S)$ as follows:
\begin{itemize}
	\item The objects of $\ho\Span^R(\Sch_S)$ are separated $S$-schemes.
	\item The set of morphisms from $X$ to $Y$ is $\pi_0 \Corr^R_S(X, Y)$.
	\item The identity morphism at $X$ is $[\Gamma^R(\id_X)] \in\pi_0 \Corr^R_S(X,X)$.
	\item The composition law is given by $\pi_0$ of the composition law~\eqref{eq:composition in Corr^R}.
\end{itemize}
It is moreover easy to show that the morphisms $\Gamma^R$ defined in \sssecref{sssec:GammaR} assemble into a functor
\[
\Gamma^R\colon \Sch_S \to \ho\Span^R(\Sch_S).
\]

For any full subcategory $C\subset \Sch_S$, we denote by $\ho\Span^R(C)$ the corresponding full subcategory of $\ho\Span^R(\Sch_S)$. 

\begin{rem}\label{rem:CorrR}
	By construction, $\ho\Span^R(C)$ is enriched in the homotopy category $\ho\Spc$. If $R$ is an $\sA_\infty$-ring spectrum, we expect that with more effort one can construct an $\infty$-category $\Span^R(C)$ with mapping spaces $\Corr^R_S(X,Y)$, whose homotopy category is $\ho\Span^R(C)$; this explains our notation for the latter category. In our two main examples, when $C$ is the category of smooth schemes over a field and $R=H\Z$ or $R=H\tilde\Z$, we will see that the spaces $\Corr^R_S(X,Y)$ are always discrete, so that $\Span^R(C) = \ho\Span^R(C)$.
\end{rem}

\sssec{}
We note that the category $\ho\Span^R(\Sch_S)$ is semiadditive, with the sum given by the disjoint union of schemes. In fact, we have canonical equivalences of spaces
\begin{gather*}
\Corr^R_S(X,\emptyset) \simeq * \simeq \Corr^R_S(\emptyset,Y),\\
\Corr^R_S(X,Y_1\coprod Y_2) \simeq \Corr^R_S(X,Y_1) \times \Corr^R_S(X, Y_2),\\
\Corr^R_S(X_1\coprod X_2,Y) \simeq \Corr^R_S(X_1,Y) \times \Corr^R_S(X_2, Y).
\end{gather*}

\begin{rem}
	In general, $\Corr^R_S(-,Y)$ is not a Nisnevich sheaf. Indeed, Calmès and Fasel show in \cite[Example 5.12]{Calmes:2014ab} that $\Corr^{\hzmw}_k(-,\A^1_k-\{0,1\})$ is not a Nisnevich sheaf.
\end{rem}

\sssec{The functor to $R$-modules} 
For $R\in \SH(S)$ a homotopy associative ring spectrum, we define a functor
\begin{equation*}\label{eqn:to-R-modules}
M^R\colon \ho\Span^R(\Sm_S) \to \Mod_R(\ho\SH(S)),\quad X\mapsto R\otimes \Sigma^\infty_\T X_+.
\end{equation*}
We shall use the fact that the four functors $f^*$, $f_*$, $f_!$, and $f^!$ preserve $R$-modules, in the sense that they lift canonically from $\ho\SH(-)$ to $\Mod_R(\ho\SH(-))$.

Let $Z\subset X\times_SY$ and $\alpha\in R^\BM(Z/X)$ define a finite $R$-correspondence from $X$ to $Y$.
 Consider the diagram
\begin{equation*}
\begin{tikzcd}
 & Z \ar[swap]{dl}{f} \ar{dr}{g} & \\
 X\ar{dr}[swap]{p}  & & Y \ar{dl}{q} \\
 & S. & 
\end{tikzcd}
\end{equation*}
As in~\sssecref{sssec:purity-transformation}, $\alpha$ gives rise to a natural transformation 
\[ 
\alpha\colon f^* \to f^!\colon \Mod_R(\ho\SH(X)) \to \Mod_R(\ho\SH(Z)).
\]
The functor $M^R$ then sends $(Z,\alpha)$ to the composition
\[
R\otimes \Sigma^\infty_\T X_+\simeq p_!p^!R \xrightarrow{\text{unit}} p_!f_*f^*p^!R \xrightarrow\alpha p_!f_!f^!p^!R \simeq q_!g_!g^!q^!R \xrightarrow{\text{counit}} q_!q^!R\simeq R\otimes \Sigma^\infty_\T Y_+.
\]
Compatibility with composition is a straightforward verification.

\begin{rem}\label{rem:CorrR2}
	If $R$ is an $\sA_\infty$-ring spectrum, then the four functors $f^*$, $f_*$, $f_!$, and $f^!$ lift canonically from $\SH(-)$ to $\Mod_R(\SH(-))$, and the above construction actually defines an $\ho\Spc$-enriched functor
\[
M^R\colon \ho\Span^R(\Sm_S) \to \ho\Mod_R(\SH(S)).
\]
	Following Remark~\ref{rem:CorrR}, we expect that it can be refined to a functor of $\infty$-categories 
	\[M^R\colon \Span^R(\Sm_S) \to \Mod_R(\SH(S)).\]
\end{rem}

\sssec{The symmetric monoidal structure}\label{sssec:CorrR-monoidal}
Let $R\in\SH(S)$ be a homotopy \emph{commutative} ring spectrum. Then the category $\ho\Span^R(\Sch_S)$ acquires a symmetric monoidal structure given on objects by $X\otimes Y=X\times_SY$. On morphisms, one uses the external pairing
\[
R^\BM(Z/X) \times R^\BM(Z'/X') \to R^\BM(Z\times_SZ'/X\times_SX').
\]
The compatibility between this pairing and composition of finite $R$-correspondences uses the commutativity of $R$. 
Furthermore, the functor $M^R\colon \ho\Span^R(\Sm_S) \to \Mod_R(\ho\SH(S))$ admits a canonical symmetric monoidal structure. We omit the somewhat tedious details.

\begin{rem}\label{rem:CorrR3}
	If $R$ is an $\sE_{n+1}$-ring spectrum ($1\leq n\leq\infty$), we expect an $\sE_n$-monoidal structure on the $\infty$-category $\Span^R(\Sm_S)$ and on the functor $M^R\colon \Span^R(\Sm_S)\to \Mod_R(\SH(S))$ (see Remark~\ref{rem:CorrR2}).
\end{rem}

\sssec{}
\label{sssec:borel-moore vs support}
 We can express finite $R$-correspondences in terms of twisted cohomology with support.
 Let $Z \subset X\times_SY $ be a reduced subscheme that is finite and universally open over $X$, and assume that $Y$ is \emph{smooth} over $S$. As explained in \sssecref{sssec:bm-support}, there is a canonical equivalence
 \begin{equation} \label{eq:bm-vs-supp}
 R^{\BM}(Z/X) \simeq R_Z(X \times_S Y, \pi_Y^* \Omega_{Y/S}),
 \end{equation}
 where $\pi_Y\colon X\times_SY\to Y$ is the projection. This equivalence is moreover natural in $Z$ by \sssecref{sssec:bm-support}(2), so that
 \[
 \Corr^R_S(X,Y) \simeq \colim_{Z\subset X\times_S Y} R_Z(X\times_S Y, \pi_Y^*\Omega_{Y/S}).
 \]

\begin{ex}
Suppose that $S$ is regular Noetherian and let $X,Y\in\Sm_S$. Then $\mathrm{KGL}^\BM(Z/X)$ is the $G$-theory space $G(Z)$,
	and hence
	\[
	\Corr^\mathrm{KGL}_S(X,Y) = \colim_{Z\subset X\times_SY} G(Z).
	\]
	Alternatively, by~\eqref{eq:bm-vs-supp} and the continuity of $K$-theory, $\Corr^\mathrm{KGL}_S(X,Y)$ is the $K$-theory space of the stable $\infty$-category of perfect complexes on $X\times_SY$ supported on a subscheme finite and equidimensional over $X$.
\end{ex}

\sssec{}
We observe that the notion of finite $R$-correspondence between smooth $S$-schemes depends only on the very effective cover of $R$ (in the sense of Spitzweck–\O stv\ae r \cite{SpitzweckOstvaer}).

\begin{prop}
	Let $S_0$ be a Noetherian scheme of finite Krull dimension with perfect residue fields and let $R\in\SH(S_0)$. Then the very effective cover $\tilde f_0 R\to R$ induces an equivalence 
	\[\Corr^{\tilde f_0R}_S(X,Y)\simeq \Corr^R_S(X,Y)\]
	for every essentially smooth $S_0$-scheme $S$ and every $X,Y\in\Sm_S$.
\end{prop}

\begin{proof}
	It will suffice to show that the map
	\[(\tilde f_0R)^\BM(Z/X) \to R^\BM(Z/X)\] 
	is an equivalence for every $X\in \Sm_S$ and every finite smoothable morphism $Z\to X$. By standard limit arguments, we can reduce to the case $S=S_0$.
	By~\eqref{purity3}, it is enough to show that the map $(\tilde f_0R)_Z(V,\xi) \to R_Z(V,\xi)$ is an equivalence for every $V\in \Sm_S$, $\xi\in K(V)$ of rank $r$, and $Z\subset V$ fiberwise of codimension $\geq r$. Since the question is local on $V$, we can assume that $\xi$ is pulled back from $S$, so that
	\[
	R_Z(V,\xi)=\Maps_{\SH(S)}(\Sigma^\infty_\T(V/V-Z), \Sigma^\xi R)
	\]
	and similarly for $\tilde f_0R$. Since $\tilde f_r \Sigma^\xi R\simeq \Sigma^\xi\tilde f_0R$, it remains to show that $\Sigma^\infty_\T(V/V-Z)$ is very $r$-effective. By \cite[Proposition B.3]{norms} and the assumptions on $S$, we may assume that $S$ is the spectrum of a perfect field. In this case, $Z$ admits a finite stratification by smooth schemes and the result is easily proved by induction using the purity isomorphism.
\end{proof}

\sssec{} 
In case $k$ is a field and $R=H\Z$ or $R=\hzmw$, the hypothetical $\infty$-category $\Span^R(\Sm_k)$ happens to be a $1$-category, i.e., it is equivalent to its homotopy category $\ho\Span^R(\Sm_k)$.
This is a special case of the following proposition. We refer to \cite[Section 3]{BachmannSlices} for the definition of the effective homotopy $t$-structure.

\begin{prop}\label{prop:discrete}
Let $k$ be a perfect field and let $R\in \SH(k)$ be a motivic spectrum in the heart of the effective homotopy $t$-structure.
For any essentially smooth $k$-scheme $S$ and $X,Y\in \Sm_S$, the $\infty$-groupoid $\Corr^R_S(X,Y)$ is discrete.
\end{prop}

\begin{proof} 
Using the description of $\Corr^R_S(X,Y)$ given in~\sssecref{sssec:borel-moore vs support}, it suffices to show that the $\infty$-groupoid
\[
R_Z(X\times_SY, \pi_Y^*\Omega_{Y/S})
\]
is discrete for any $Z\subset X\times_SY$ finite over $X$.
By a standard limit argument, we can assume that $S$ is smooth over $k$. The result then follows from Lemma~\ref{lem:discrete} below.
\end{proof}

\begin{lem}\label{lem:discrete}
	Let $k$ be a perfect field, $V$ a smooth $k$-scheme, $\xi\in K(V)$, and $Z\subset V$ a closed subscheme of codimension $\geq \rk\xi$.
	Let $R\in \SH(k)$ be a motivic spectrum in the heart of the effective homotopy $t$-structure. Then the $\infty$-groupoid $R_Z(V,\xi)$ is discrete.
\end{lem}

\begin{proof}
	Suppose first that $Z$ is smooth.
	We then have the purity equivalence~\eqref{eqn:purity-support}
	\[
	R_Z(V,\xi) \simeq R(Z,\xi-\sN_{Z/X}).
	\]
	Let $\zeta = \xi-\sN_{Z/X}$. By the assumption on $\xi$, we have $\rk\zeta\leq 0$. 
	The assumption on $R$ means that $R$ is right orthogonal to $\SH^\eff_{\geq 1}(k)$. 
	Since $Z$ is smooth, $R_Z\in \SH(Z)$ is right orthogonal to $\SH^\eff_{\geq 1}(Z)$, hence so is $\Sigma^{\zeta}R_Z$, because  $\Sigma^{-\zeta}$ is a right $t$-exact endomorphism of $\SH^\eff(Z)$. It follows at once that $R(Z,\zeta)=\Maps_{\SH(Z)}(\1_Z, \Sigma^\zeta R_Z)$ is discrete.

	If $Z$ is an arbitrary closed subscheme, we can assume that it is reduced since cohomology with support only depends on $Z_{\red}$. We will prove the claim by induction on the dimension of $Z$. If $Z$ is empty, then $R_Z(V,\xi)$ is contractible. Otherwise, since $k$ is perfect, $Z$ is generically smooth, so there is a reduced closed subscheme $Z_1 \subset Z$ of strictly smaller dimension such that $ Z-Z_1$ is smooth. By~\eqref{eqn:E_Z(X)}, we have a fiber sequence (of grouplike $\Einfty$-spaces)
	 \[
	 R_{Z_1}(V, \xi) \to R_Z(V, \xi) \to R_{Z-Z_1}(V-Z_1, \xi).
	 \]
	 By the induction hypothesis, $R_{Z_1}(V, \xi)$ is discrete since $\rk \xi \leq \codim(Z,V)\leq \codim(Z_1,V)$.
	 Since $Z-Z_1$ is smooth, we have already proven that $R_{Z-Z_1}(V-Z_1, \xi)$ is discrete. It follows that $R_Z(V,\xi)$ is also discrete.
\end{proof}

\begin{rem}\label{rem:discrete}
	The converse of Proposition~\ref{prop:discrete} also holds trivially: if $\Corr^R_k(X,Y)$ is discrete for all $X,Y\in \Sm_k$, then in particular $R(X)\simeq \Corr^R_k(X,\Spec k)$ is discrete for all $X\in \Sm_k$, i.e., $\tilde f_0R$ belongs to the heart of the effective homotopy $t$-structure. 
	Thus, the hypothetical $\infty$-category $\Span^R(\Sm_k)$ is a $1$-category if and only if $\tilde f_0R\in\SH^\eff(k)^\heartsuit$, in which case $R$ is necessarily an algebra over $\underline\pi{}_0^\eff(\1)\simeq\hzmw$. In particular, for more general $R$, there is no hope to recover $R$-modules from the $1$-category $\ho\Span^R(\Sm_k)$. 
\end{rem}

\sssec{}
\label{sssec:composition with support}
It will be useful to have a description of the composition in $\ho\Span^R(\Sm_S)$ in terms of cohomology with support.
Suppose that $X, Y, T$ are smooth $S$-schemes and $Z \subset X \times_S Y$ and $Z' \subset Y \times_S T$ are reduced subschemes that are finite and universally open over $X$ and $Y$ respectively. We will refer to the diagram
\begin{equation*} \label{xyz}
\begin{tikzcd}
X \times_S T \ar{dd}{q_X} \ar{rr}{q_T} & & T \\
 & X \times_S Y \times_S T \ar[swap, "p_{XT}"]{ul} \ar{d}{p_{XY}} \ar{r}{p_{YT}} & Y \times_S T \ar[swap, "s_Y"]{d} \ar{u}{s_T} \\
X & X \times_S Y \ar{r}{r_Y} \ar[swap, "r_X"]{l} & Y.
\end{tikzcd}
\end{equation*} 
 Recall from \sssecref{sssec:CorrR-composition} that $Z' \circ Z \subset X \times_S T$ is the reduced subscheme $p_{XT}(Z'')\subset X\times_ST$ where $Z''=Z\times_YZ'\subset X\times_SY\times_ST$. 
We define the pairing
\begin{equation*} \label{eq:comp-w-supports}
\theta\colon R_Z(X \times_S Y, r_Y^*\Omega_{Y/S}) \times  R_{Z'}(Y \times_S T, s_T^*\Omega_{T/S}) \rightarrow  R_{Z' \circ Z}(X \times_S T, q_T^*\Omega_{T/S})
\end{equation*}
as the composition
\begin{equation*}\label{eqn:composition}
\begin{tikzcd}
 R_Z(X \times_S Y, r_Y^*\Omega_{Y/S}) \times  R_{Z'}(Y \times_S T, s_T^*\Omega_{T/S}) \ar{d}{p_{XY}^* \times p_{YT}^*} \\
 R_{Z\times_ST}(X \times_S Y \times_S T, p_{XY}^*r_Y^*\Omega_{Y/S}) \times  R_{X\times_SZ'}(X\times_S Y \times_S T, p_{YT}^*s_T^*\Omega_{T/S}) \ar{d}{\mu} \\
 R_{Z''}(X \times_S Y \times_S T, p_{XY}^*r_Y^*\Omega_{Y/S} + p_{YT}^*s_T^*\Omega_{T/S}) \ar{d}{p_{XT!}} \\
 R_{Z' \circ Z}(X \times_S T,  q_T^*\Omega_{T/S}),
\end{tikzcd}
\end{equation*}
where $\mu$ is the cup product and $p_{XT!}$ is the Gysin map~\eqref{eq:gysin-in-coh}. 
More succinctly,
\[
\theta(x,y)= p_{XT!}(p_{XY}^*x \cup p_{YT}^*y).
\]
We have the following comparison with the pairing $\theta^\BM$ defined in~\sssecref{sssec:CorrR-composition}:

\begin{lem}\label{lem:compositions-agree}
	The following diagram commutes, where the horizontal equivalences are instances of \eqref{eq:bm-vs-supp}:
	\begin{equation*}
	\begin{tikzcd}
	 R_Z(X \times Y, \Omega_{Y}) \times  R_{Z'}(Y \times T, \Omega_{T}) \ar{r}{\simeq} \ar{d}[swap]{\theta} & R^{\BM}(Z/X) \times R^{\BM}(Z'/Y) \ar{d}{\theta^\BM}\\
	 R_{Z \circ Z'}(X \times T, \Omega_{T}) \ar{r}{\simeq} & R^{\BM}(Z \circ Z'/X).
	\end{tikzcd}
	\end{equation*}
\end{lem}

\begin{proof} 
Recall from \sssecref{sssec:bm-support} that we can factor the cup product as follows:
\[
\begin{tikzcd}
	R_{Z\times T}(X \times Y \times T, \Omega_{Y}) \times  R_{X\times Z'}(X\times Y \times T, \Omega_{T}) \ar{d}[swap]{\id\times i^*} \ar{dr}{\mu} & \\
	R_{Z\times T}(X \times Y \times T, \Omega_{Y}) \times  R_{Z''}(Z \times T, \Omega_{T}) \ar{r}{\bar\mu} & 
	R_{Z''}(X \times Y \times T, \Omega_{Y\times T}),
\end{tikzcd}
\]
where $i\colon Z\times T\hook X\times Y\times T$.
In the following diagram, the left column is $\theta$ and the right column is $\theta^\BM$:
\[
\begin{tikzcd}
R_Z(X \times Y, \Omega_{Y}) \times R_{Z'}(Y \times T, \Omega_{T}) \ar{r}{\simeq} \ar{d}[swap]{\id \times i^*p_{YT}^*} & R^{\BM}(Z/X) \times R^{\BM}(Z'/Y) \ar{d}{\id\times g^*} \\
R_Z(X \times Y, \Omega_{Y}) \times R_{Z''}(Z \times T, \Omega_{T}) \ar{r}{\simeq} \ar{d}[swap]{p_{XY}^*\times \id}& R^{\BM}(Z/X) \times R^{\BM}(Z''/Z) \ar{dd}{\mu^\BM}\\
R_{Z\times T}(X \times Y \times T, \Omega_{Y}) \times R_{Z''}(Z \times T, \Omega_{T}) \ar{d}[swap]{\bar\mu} & \\
R_{Z''}(X \times Y \times T, \Omega_{Y\times T}) \ar{r}{\simeq} \ar{d}[swap]{p_{XT!}} & R^{\BM}(Z''/X) \ar{d}{p_*}\\
R_{Z \circ Z'}(X \times T, \Omega_{T}) \ar{r}{\simeq} & R^{\BM}(Z\circ Z'/X).
\end{tikzcd}
\]
The three rectangles commute by \sssecref{sssec:bm-support} (1), (3), and (2), respectively.
\end{proof}


\ssec{From framed correspondences to finite $R$-correspondences}
\label{ssec:functor}

Let $R$ be a homotopy associative ring spectrum. We will construct a canonical functor
$$\Phi^R \colon \ho\Span^\fr(\Sch_S) \longrightarrow \ho\Span^R(\Sch_S),$$
where $\Span^\fr(\Sch_S)$ is the $\infty$-category of framed correspondences constructed in \cite{EHKSY1}.

\sssec{}
For $S$-schemes $X$ and $Y$, we define a map
\begin{equation}\label{eqn:Phi-on-morphisms}
\Phi^R\colon \Corr^\fr_S(X,Y) \to \Corr^R_S(X,Y)
\end{equation}
as follows.
A framed correspondence from $X$ to $Y$
is given by a span
\begin{equation*}
\begin{tikzcd}
 & Z \ar[swap]{dl}{f} \ar{dr}{g} & \\
 X  & & Y
\end{tikzcd}
\end{equation*}
where $f$ is finite syntomic, together with a trivialization $\tau \in \Maps_{K(Z)}(0, \sL_f)$.

Since the morphism $f$ is finite syntomic, it has a fundamental class
\[
\eta_f \in \1^{\BM}(Z/X,\sL_f)
\]
(defined using the canonical factorization  $Z\hook \bV(f_*\sO_Z)\to X$, see \sssecref{sssec:fund-tr}). 
We will also denote by 
\[
\eta_f \in R^{\BM}(Z/X,\sL_f)
\]
its image by the map $\1^{\BM}(Z/X,\sL_f)\to R^{\BM}(Z/X,\sL_f)$ induced by the unit $\1_S\to R$.
Applying the trivialization $\tau$, we get an element $\tau_*(\eta_f) \in R^{\BM}(Z/X)$. 

The map $(f, g) \colon Z \to X \times_S Y $ is finite; we denote by $V\subset X \times_S Y$ its reduced image. Note that $V$ is finite and universally open over $X$. Using the proper pushforward in Borel–Moore homology, we obtain  
$(f, g) _*(\tau_*(\eta_f)) \in R^{\BM}(V/X).$ 
 This construction defines a map
 \[
 \Maps_{K(Z)}(0,\sL_f) \to R^\BM(V/X), \quad \tau \mapsto (f,g)_*(\tau_*(\eta_f)).
 \]
 Taking the colimit over the groupoid of finite syntomic spans from $X$ to $Y$, we obtain the desired map~\eqref{eqn:Phi-on-morphisms}.

\begin{prop}\label{prop:PhiR}
The maps~\eqref{eqn:Phi-on-morphisms} define an $\ho\Spc$-enriched functor 
\[ \Phi^R \colon \ho\Span^\fr(\Sch_S) \longrightarrow \ho\Span^R(\Sch_S) \]
such that the following triangle commutes:
\[
\begin{tikzcd}
	\Sch_S \ar{r}{\Gamma^R} \ar{d}[swap]{\gamma} & \ho\Span^R(\Sch_S). \\
	\ho\Span^\fr(\Sch_S) \ar{ur}[swap]{\Phi^R} &
\end{tikzcd}
\]
\end{prop}

\begin{proof}
It is clear that $\Phi^R(\gamma(f))=\Gamma^R(f)$ for any $S$-morphism $f$, and in particular $\Phi^R$ preserves identity morphisms.
Let $\alpha = (Z, f, g, \tau) \in \Corr^\fr_S(X, Y)$ and $\beta = (Z', h, s, \tau') \in \Corr^\fr_S(Y, T)$ be framed correspondences, where  $\tau \in \Maps_{K(Z)}(0, \sL_f)$ and $\tau' \in \Maps_{K(Z')}(0, \sL_h)$, and form the composite $2$-span
 \begin{equation*}
 \begin{tikzcd}[column sep={4.5em,between origins}]
 & & Z'' \ar{dr}{g'} \ar{dl}[swap]{h'} &  & \\
 & Z \ar{dr}{g} \ar{dl}[swap]{f} &  &  Z' \ar{dl}[swap]{h} \ar{dr}{k}  & \\
 X & & Y & & T\rlap.
 \end{tikzcd}
 \end{equation*}
 The composition $\beta\circ\alpha$ is then given by 
 $$\beta\circ\alpha = (Z'', f \circ h', k \circ g', \sigma) \in \Corr^\fr_S(X, T),$$
  where $\sigma\in\Maps_{K(Z'')}(0,\sL_{f\circ h'})$ is the composite
 $$0 \xrightarrow{\tau\oplus\tau'} h'^*\sL_f \oplus g'^*\sL_{h} \simeq \sL_{f\circ h'}.$$
 
 We want to compare $\Phi^R(\beta) \circ \Phi^R(\alpha)$ and $\Phi^R(\beta \circ \alpha)$.
 We first note the following equations between fundamental classes:
 \[
 g^*(\tau'_*\eta_{h})\cdot \tau_*\eta_f  = 
  g'^* (\tau')_* \eta_{h'}\cdot \tau_*\eta_f = \sigma_*( \eta_{h'}\cdot\eta_f ) = \sigma_* \eta_{f\circ h'}.
 \]
 Here the first equality is the stability of fundamental classes under tor-independent base change, the second holds by definition of $\sigma$, and the last is the associativity of fundamental classes \cite[Definition~2.3.6]{DJK}.
Let $V\subset X\times_SY$ be the image of $Z$ and $V'\subset Y\times_ST$ the image of $Z'$.
We now consider the following diagram in which all parallelograms are Cartesian:
\[
\begin{tikzcd}[row sep=1em, column sep=3em]
	&& Z'' \ar{d}[swap]{r} \ar{ddl} \ar{ddr} && \\
	&& W \ar{dl} \ar{d}[swap]{q} \ar{ddr} && \\
	&Z \ar{d}[swap]{p} \ar{ddr}[description]{g} \ar{ddl} & V'' \ar{dl} \ar{dr} & Z' \ar{d}[swap,xshift=2pt]{p'} \ar{ddl} \ar{ddr} & \\
	&V \ar{dr}[swap]{v} \ar{dl} &  & V' \ar{dl} \ar{dr} & \\
	X && Y && T\rlap.
\end{tikzcd}
\]
For any $z\in R^\BM(Z/X)$ and $z'\in R^\BM(Z'/Y)$, we have the following equivalences in $R^\BM(V''/X)$, where the parenthetical justifications refer to \cite[1.2.8]{deglise2017bivariant}:
\begin{align*}
v^*p'_*(z')\cdot p_*(z) & = q_*(p^*v^*p'_*(z')\cdot z) & \text{(projection formula)} \\
&=q_*(g^*p'_*(z')\cdot z) & \text{(composition)}\\
&=q_*(r_*g^*(z')\cdot z)& \text{(base change)} \\
&= q_* r_* (g^*(z')\cdot z)& \text{(compatibility with pushforwards)} \\
&= (q\circ r)_*(g^*(z')\cdot z). & \text{(composition)}
\end{align*}
Plugging in $z=\tau_*\eta_f$ and $z'=\tau'_*\eta_h$ and pushing forward the result to $R^\BM(V'\circ V/X)$ gives the desired equivalence
\[
\Phi^R(\beta) \circ \Phi^R(\alpha)\simeq \Phi^R(\beta \circ \alpha).
\]
To see that $\Phi^R$ is indeed an $\ho\Spc$-enriched functor, we must show that this equivalence is natural in the pair $(\alpha,\beta)\in \Corr^\fr_S(X, Y)\times \Corr^\fr_S(Y, Z)$. This is essentially obvious from the construction, using the functoriality of fundamental classes discussed in \sssecref{sssec:eta-functoriality} and \sssecref{sssec:eta-composition}.
\end{proof}

The following corollary is a variant of \cite[Theorem 10.1]{DruzhininKolderup}.

\begin{cor}\label{cor:strict-A1-invariance}
	Let $k$ be a perfect field and $R\in\SH(k)$ a homotopy associative ring spectrum. Let $\sF$ be an $\A^1$-invariant presheaf of abelian groups on $\ho\Span^R(\Sm_k)$ that preserves finite products. Then the Nisnevich sheaf $L_\nis \sF$ on $\Sm_k$ is strictly $\A^1$-invariant and the canonical map $H^i_\zar(-,L_\zar\sF)\to H^i_\nis(-,L_\nis\sF)$ is an isomorphism for all $i\geq 0$.
\end{cor}

\begin{proof}
	This follows from the existence of the functor $\Phi^R\colon \Span^\fr(\Sm_k)\to \ho\Span^R(\Sm_k)$ and \cite[Theorem 3.4.11]{EHKSY1}.
\end{proof}

\sssec{}
Suppose that $R$ is homotopy commutative. As explained in \sssecref{sssec:CorrR-monoidal}, $\ho\Span^R(\Sm_S)$ is then a symmetric monoidal category. Recall that $\Span^\fr(\Sm_S)$ is also a symmetric monoidal $\infty$-category. One can easily check that the functor
\[
\Phi^R\colon \ho\Span^\fr(\Sm_S) \to \ho\Span^R(\Sm_S)
\]
can be uniquely promoted to a symmetric monoidal functor in such a way that the natural equivalence $\Phi^R\circ \gamma \simeq \Gamma^R$ is monoidal.

\sssec{}
We now relate the functor $\Phi^R$ to the free $R$-module functor:

\begin{prop}\label{prop:DM-vs-Mod}
	Let $R\in\SH(S)$ be a homotopy associative ring spectrum.
	 Then the following diagram of $\ho\Spc$-enriched categories commutes:
	\[
	\begin{tikzcd}
	\ho\Span^\fr(\Sm_S) \ar{r}{\gamma_*\Sigma^\infty_{\T,\fr}} \ar{d}{\Phi^R} & \ho\SH(S) \ar{d}{ R\otimes -} \\
	\ho\Span^R(\Sm_S) \ar{r}{M^R} & \Mod_R(\ho\SH(S)).
	\end{tikzcd}
	\]
	Furthermore, if $R$ is homotopy commutative, this square commutes in the $2$-category of symmetric monoidal $\ho\Spc$-enriched categories.
\end{prop}

\begin{proof}
By definition of these functors, we have given isomorphisms $\M^R\Phi^R(X) \simeq R\otimes \Sigma^\infty_\T X_+\simeq R\otimes \gamma_*\Sigma^\infty_{\T,\fr}X$ for every $X\in \Sm_S$.
Moreover, when $R$ is homotopy commutative, these isomorphisms trivially intertwine the monoidal structures of these functors.
It thus remains to show that the following square commutes for every $X,Y\in \Sm_S$:
\begin{equation}\label{eqn:R-module-square}
	\begin{tikzcd}
		\Corr^\fr_S(X,Y) \ar{r} \ar{d} & \Maps_{\SH(S)}(\Sigma^\infty_\T X_+, \Sigma^\infty_\T Y_+) \ar{d} \\
		\Corr^R_S(X,Y) \ar{r} & \Maps_{\SH(S)}(R\otimes \Sigma^\infty_\T X_+, R\otimes \Sigma^\infty_\T Y_+).
	\end{tikzcd}
\end{equation}
By Theorem~\ref{thm:fund-vs-recog}, the top horizontal map in~\eqref{eqn:R-module-square} is the fundamental transfer $\tr_\eta$.
Let $\phi=(Z,f,g,\tau)$ be a framed correspondence from $X$ to $Y$, and let $V\subset X\times_SY$ be the reduced image of $Z$:
\begin{equation*}
\begin{tikzcd}
 & Z \ar[swap]{dl}{f} \ar{dr}{g} \ar{d}{r} & \\
 X\ar{dr}[swap]{p}  & V \ar{l}{u} \ar{r}[swap]{v} & Y \ar{dl}{q} \\
 & S. & 
\end{tikzcd}
\end{equation*}
Then $\Phi^R(\phi)=(V,\alpha)$ for some $\alpha\in R^\BM(V/X)$, inducing a natural transformation $\alpha\colon u^* \to u^!$ in $R$-modules. In the following diagram, the top row is $\tr_\eta(\phi)$ (in the form described in Remark~\ref{rem:homological-gysin}), while the bottom row is $M^R\Phi^R(\phi)$:
\[
\begin{tikzcd}
p_!p^!\1_S \ar{d} \ar{r}{\text{unit}} & p_!f_*f^*p^!\1_S \ar{r}{\tau\circ\pur_f} \ar{d} & p_!f_!f^!p^!\1_S \ar{r}{\simeq} \ar{d} & q_!g_!g^!q^!\1_S \ar{r}{\text{counit}} \ar{d} & q_!q^!\1_S \ar{d} \\
p_!p^!R \ar{r}{\text{unit}} \ar[equal]{d} & p_!f_*f^*p^!R \ar{r}{\tau\circ\pur_f} & p_!f_!f^!p^!R \ar{d}{\text{counit}} \ar{r}{\simeq} & p_!g_!g^!p^!R \ar{d}{\text{counit}} \ar{r}{\text{counit}} & q_!q^!R \ar[equal]{d} \\
p_!p^!R \ar{r}{\text{unit}} & p_!u_*u^*p^!R \ar{r}{\alpha} \ar{u}{\text{unit}} & p_!u_!u^!p^!R \ar{r}{\simeq} & p_!v_!v^!p^!R \ar{r}{\text{counit}} & q_!q^!R.
\end{tikzcd}
\]
The square involving $\alpha$ commutes by definition of $\alpha$. The commutativity of the boundary of this diagram witnesses the commutativity of the square~\eqref{eqn:R-module-square}.
\end{proof}

\begin{rem}
	If $R$ is $\sA_\infty$, we can replace the lower right corner in Proposition~\ref{prop:DM-vs-Mod} by $\ho\Mod_R(\SH(S))$.
	 Continuing Remarks~\ref{rem:CorrR}, \ref{rem:CorrR2}, and \ref{rem:CorrR3}, we moreover expect that this square can be promoted to a commuting square of $\infty$-categories, and of $\sE_n$-monoidal $\infty$-categories if $R$ is $\sE_{n+1}$.
\end{rem}


\ssec{Voevodsky correspondences and Milnor–Witt correpondences}
\label{ssec:VMW}

We show that the $\infty$-category of finite $H\Z$-correpondences (resp.\ of finite $\hzmw$-correspondences) recovers Voevodsky's category of finite correspondences \cite[Lecture 1]{mvw} (resp.\ Calmès and Fasel's category of finite Milnor–Witt correspondences \cite{Calmes:2014ab}). We then show that the functor $\Phi^{H\Z}$ (resp.\ $\Phi^{\hzmw}$) recovers the functor $\cyc$ constructed in \cite[\sectsign 5.3]{EHKSY1} (resp.\ the functor of Déglise and Fasel \cite[Proposition 2.1.12]{DegliseFasel}).

\sssec{Reminders on motivic cohomology}
Over a Dedekind domain $D$, we will consider the motivic cohomology spectrum $H\Z\in \SH(D)$ constructed by Spitzweck \cite[Definition 4.27]{SpitzweckHZ}.
It is an oriented $\Einfty$-ring spectrum that represents Bloch–Levine motivic cohomology. In particular, for an essentially smooth $D$-scheme $X$ and $\xi\in K(X)$ of rank $r$, we have
\begin{equation}\label{eqn:Voevodsky-vs-Bloch}
H\Z(X,\xi) \simeq z^r_{\zar}(X,*),
\end{equation}
where $z_\zar^r(X,*)$ denotes the sheafification of Bloch's cycle complex $z^r(X,*)$ with respect to the Zariski topology on $\Spec D$. This identification is natural in $X$, where the functoriality of Bloch's cycle complex comes from Levine's moving lemma~\cite{LevineChow}.
If $Z\subset X$ is a closed subscheme of codimension $c$, the localization theorem \cite[Theorem 1.7]{LevineLocalization} implies that
\begin{equation}\label{eqn:HZ-support}
H\Z_Z(X,\xi) \simeq z^{r-c}_{\zar}(Z,*).
\end{equation}

Recall that $H\Z$ belongs to the heart of the effective homotopy $t$-structure on $\SH(D)$ \cite[Lemma 13.6]{norms}.
Being the zeroth slice of the sphere spectrum \cite[Theorem B.4]{norms}, the spectrum $H\Z$ admits in fact a unique $\Einfty$-ring structure with given unit. 

When $D$ is a field, $H\Z$ coincides with Voevodsky's motivic cohomology spectrum, but this is not known in general.
In this case, Bloch's cycle complex admits an $\Einfty$-ring structure compatible with the intersection of cycles \cite[\sectsign 5]{Bloch}, and the equivalence~\eqref{eqn:Voevodsky-vs-Bloch} is multiplicative.

\sssec{$H\Z$-correspondences vs.\ Voevodsky correspondences} We let $\Cor_S$ denote Voevodsky's category of finite correspondences between smooth separated $S$-schemes, as defined in \cite[Appendix 1A]{mvw}.

\begin{lem}\label{lem:gysin-HZ}
	Let $S$ be the spectrum of a Dedekind domain, $f\colon X\to Y$ a morphism between essentially smooth $S$-schemes, $Z\subset X$ a closed subscheme flat over $S$ such that the restriction of $f$ to $Z$ is finite, and $\xi\in K(Y)$ of rank $\codim(f(Z),Y)$. Then the Gysin map
	\[
	f_!\colon H\Z_Z(X, f^*\xi+\sL_f) \to H\Z_{f(Z)}(Y,\xi)
	\]
	agrees with the pushforward of codimension $0$ cycles $f_*\colon z^0(Z) \to z^0(f(Z))$ under the identification~\eqref{eqn:HZ-support}.
\end{lem}

\begin{proof}
	Let $\eta\in S$ be the generic point. Since $f(Z)$ is flat over $S$, the pullback $z^0(f(Z))\to z^0(f(Z)\times_S\eta)$ is an isomorphism, so we may assume that $S$ is the spectrum of a field $k$.
		By limit arguments, we can assume $k$ perfect and $Y$ smooth over $k$. Replacing $Y$ by an open subscheme, we can further assume that $Z$ and $f(Z)$ are smooth over $k$. Since the Gysin map is compatible with purity isomorphisms, we are reduced to the following claim: if $L/K$ is a finite extension of finitely generated fields over $k$, the Gysin map $\Z\simeq H\Z(\Spec L,\sL_{L/K}) \to H\Z(\Spec K)\simeq \Z$ is multiplication by $[L:K]$. This is a special case of \cite[Example 3.2.9(1)]{deglise2011orientation}.
\end{proof}

\begin{prop}\label{prop:CorrHZ-comparison}
	Let $S$ be essentially smooth over a Dedekind domain. Then the symmetric monoidal $\infty$-category $\Span^{H\Z}(\Sm_S)$ is a $1$-category and is equivalent to $\Cor_S$.
\end{prop}

\begin{proof}
	For smooth $S$-schemes $X$ and $Y$, we have
		$$
		\Corr^{H\Z}_S(X, Y) \simeq  \colim_{Z \subset X \times_S Y} H\Z_Z(X \times_S Y, \pi_Y^* \Omega_{Y/S}),
		$$
	by~\sssecref{sssec:borel-moore vs support}.
	It follows from~\eqref{eqn:HZ-support} that 
	\[
	H\Z_Z(X \times_S Y, \pi_Y^* \Omega_{Y/S})\simeq z^0(Z,*) \simeq \bigoplus_{Z^{(0)}}\Z.
	\] 
	Hence,
	\[
	\Corr^{H\Z}_S(X,Y) \simeq \bigoplus_{Z\subset X\times_SY} \Z,
	\]
	where the sum is taken over all integral closed subschemes of $X\times_SY$ that are finite and surjective over a component of $X$. In particular, $\Span^{H\Z}(\Sm_S)$ is a $1$-category, and its mapping spaces are the same as in Voevodsky's category.
	
	To compare the composition laws, we use the description of the composition in $\Span^{H\Z}(\Sm_S)$ via the pairing $\theta$ (Lemma~\ref{lem:compositions-agree}). The composition in $\Cor_S$ is defined in exactly the same way, except that it uses the intersection product and the pushforward of cycles instead of the cup product and the Gysin map in $H\Z$-cohomology.
	We must therefore show that these constuctions yield cycles with the same multiplicities. Since the generic points of the cycles involved lie over generic points of $S$, we can replace $S$ by its generic points and hence assume that $S$ is a field. In this case, the intersection product and the cup product agree because the isomorphism~\eqref{eqn:Voevodsky-vs-Bloch} is compatible with the multiplicative structures. The fact that the pushforwards agree is a special case of Lemma~\ref{lem:gysin-HZ}.
	Finally, the fact that the symmetric monoidal structures agree also follows from the multiplicativity of the isomorphism~\eqref{eqn:Voevodsky-vs-Bloch}.
\end{proof}

\sssec{} In \cite[\sectsign5.3]{EHKSY1}, we defined a symmetric monoidal functor
\[
\cyc\colon \Span^\fr(\Sm_S) \to \Cor_S
\]
sending a framed correspondence
\begin{equation*}
  \begin{tikzcd}
    & Z \ar[swap]{ld}{f}\ar{rd}{g}
  \\
  X & & Y
  \end{tikzcd}
\end{equation*}
to the cycle $(f,g)_*[Z]$ on $X\times_SY$, where $[Z]\in z^0(Z)$ is the fundamental cycle of $Z$. By Proposition~\ref{prop:CorrHZ-comparison}, we also have the symmetric monoidal functor
\[
\Phi^{H\Z}\colon \Span^\fr(\Sm_S) \to \Span^{H\Z}(\Sm_S)\simeq \Cor_S
\]
defined in \ssecref{ssec:functor}.

\begin{prop}\label{prop:PhiHZ}
	For $S$ essentially smooth over a Dedekind domain, there is an isomorphism of symmetric monoidal functors
	\[
	\Phi^{H\Z}\simeq \cyc\colon \Span^\fr(\Sm_S) \to \Cor_S.
	\]
\end{prop}

\begin{proof}
Note that $\Phi^{H\Z}$ and $\cyc$ send a framed correspondence to finite correspondences with the same support, so it suffices to compare their multiplicities. 
Since the generic points of their support lie over generic points of $S$ and both functors are natural in $S$, this can be done assuming that $S=\Spec k$ for some field $k$, which can moreover be assumed perfect by passing to its perfection.
	In this situation, we prove the following more general uniqueness statement: if
	\[
	\phi_1,\phi_2\colon \Span^\fr(\Sm_k) \to \Span^{H\Z}(\Sm_k)
	\]
	are symmetric monoidal functors
	that satisfy $\phi_1|\Sm_k\simeq \Gamma^{H\Z}\simeq \phi_2|\Sm_k$ and send every framed correspondence $(Z,f,g,\tau)$ to a finite correspondence with support $(f,g)(Z)$, then $\phi_1\simeq \phi_2$.
We have induced symmetric monoidal functors
	\[
	\phi_1^*,\phi_2^*\colon \SH^\fr(k) \to \DM(k)
	\]
	such that $\phi_1^*|\SH(k)\simeq \phi_2^*|\SH(k)$. By the reconstruction theorem \cite[Theorem 3.5.12]{EHKSY1} it follows that $\phi_1^*\simeq \phi_2^*$.
	To check that $\phi_1\simeq \phi_2$, it suffices to compare their effect on a framed correspondence $\alpha\in \Corr^\fr_k(\eta,Y)$ with connected support, where $\eta$ is the generic point of a smooth $k$-scheme. 
	Since $\phi_1(\alpha)$ and $\phi_2(\alpha)$ are supported on a single point, their equality can be checked modulo rational equivalence, i.e., in $\ho\DM(k)$, so we are done.
\end{proof}

\sssec{Reminders on Milnor–Witt motivic cohomology}
Over a field $k$, we will consider the Milnor–Witt motivic cohomology spectrum $\hzmw\in\SH(k)$. We adopt the definition
\[
\hzmw = \underline\pi^\eff_0(\1),
\]
where $\underline\pi^\eff_*$ are the homotopy groups in the effective homotopy $t$-structure. 
This definition is due to Bachmann \cite{BachmannSlices}, and it is known to be equivalent to that of Calmès and Fasel when $k$ is infinite perfect of characteristic not $2$ \cite{BachmannFasel}.
By definition, $\hzmw$ is an $\Einfty$-ring spectrum in the heart of the effective homotopy $t$-structure. Moreover, since the unit map $\1\to \mathrm{MSL}$ is a $\underline\pi^\eff_0$-isomorphism \cite[Example 16.34]{norms}, $\hzmw$ is uniquely $\SL$-oriented. In particular, we have \emph{Thom isomorphisms}
\[
\Sigma^\xi\hzmw_X \simeq \Sigma^{\rk\xi}_\T \Sigma^{\det\xi-\sO} \hzmw_X
\]
for any $X\in\Sm_X$ and $\xi\in K(X)$ (apply \cite[Example 16.29]{norms} to $\xi-\det\xi$).

Since the effective cover functor $f_0\colon \SH(k)\to \SH^\eff(k)$ is $t$-exact for the respective homotopy $t$-structures \cite[Proposition 4(3)]{BachmannSlices}, we have
\[
 \hzmw \simeq f_0\underline \pi_0(\1)_*.
\]
Recall Morel's computation $\underline\pi_0(\1)_* \simeq \underline K^{MW}_*$ \cite{morel-pi0}. More generally, for $X$ a smooth $k$-scheme and $\sL$ an invertible sheaf on $X$, we have an equivalence in $\SH(X)^\heartsuit$
\[
\Sigma^{\sL-\sO}\underline\pi_0(\1_X)_* \simeq \underline K^{MW}_*(\sL) = \underline K^{MW}_* \otimes_{\Z[\sO^\times]} \Z[\sL^\times]
\]
by \cite[Lemma 2.9]{Ananyevskiy}. Therefore the canonical map $\hzmw \to \underline K^{MW}_*$ and the $\SL$-orientation of $\hzmw$ induce maps of abelian groups
\begin{equation}\label{eqn:hzmw-to-kmw}
\pi_i\hzmw(X,\xi) \to H^{n-i}_\nis(X, \underline K^{MW}_n(\det\xi)),
\end{equation}
natural in $X\in \Sm_k$ and $\xi\in K(X)$, where $n=\rk\xi$.
If $k$ is perfect, we can analyze the Postnikov filtration of $f_0\underline K^{MW}_*$ using Rost–Schmid complexes (see \sssecref{sssec:rost-schmid} below), and we easily deduce that~\eqref{eqn:hzmw-to-kmw} is an isomorphism for $i=0,1$.
In particular, if $Z\subset X$ is a closed subset, then
\begin{equation}\label{eqn:hzmw-chw}
\pi_0\hzmw_Z(X,\xi) \simeq \Chw^{n}_{Z}(X,\det\xi)= H^{n}_{\nis,Z}(X, \underline K^{MW}_n(\det\xi)).
\end{equation}
By continuity, these computations remain valid over arbitrary fields. Moreover, the isomorphism~\eqref{eqn:hzmw-chw} is compatible with the multiplicative structures, since the product in Milnor–Witt K-theory (which induces the intersection product on Chow–Witt groups) is induced by the ring structure of the sphere spectrum.

\sssec{Rost–Schmid complexes}\label{sssec:rost-schmid}
Let $k$ be a perfect field, $X$ a smooth $k$-scheme, $\sL$ an invertible sheaf on $X$, and $Z\subset X$ a closed subset.
The Nisnevich cohomology of $X$ with coefficients in the sheaf $\underline K^{MW}_n(\sL)$ and with support in $Z$ can be computed using the Rost–Schmid complex $C^*_Z(X,\underline K^{MW}_n(\sL))$, given in degree $i$ by
\[
C^i_Z(X,\underline K^{MW}_n(\sL)) = \bigoplus_{x\in X^{(i)}\cap Z} K^{MW}_{n-i}(\kappa(x),\sL\otimes\omega_x),
\]
where $\omega_{x}=\omega_{\kappa(x)/\sO_{X,x}}$ \cite[Chapter 5]{Morel}. 
In particular, 
\begin{equation}\label{eqn:chw-rs}
	\Chw^n_Z(X,\sL)\simeq H^n(C^*_Z(X,\underline K^{MW}_n(\sL))).
\end{equation}

The Rost–Schmid complex is functorial for \emph{flat} morphisms in $\Sm_k$: if $f\colon Y\to X$ is flat, there is an induced map of complexes
\[
f^*\colon C^*_Z(X,\underline K^{MW}_n(\sL)) \to C^*_{f^{-1}(Z)}(Y,\underline K^{MW}_n(f^*\sL)),
\]
defined in~\cite[Corollaire~10.4.3]{FaselChowWitt}.
On the other hand, for any $f$, there is a sheaf-theoretic pullback
$$f^* \colon \Chw^n_Z(X,\sL) \to \Chw^n_{f^{-1}(Z)}(Y,f^*\sL),$$
which agrees with the pullback in $\hzmw$-cohomology (by the naturality of~\eqref{eqn:hzmw-chw}).

\begin{lem}
\label{lem:flat pullback}
The isomorphism~\eqref{eqn:chw-rs} is natural with respect to flat morphisms in $\Sm_k$. 
\end{lem}

\begin{proof}
	It suffices to show that the canonical inclusion
	\[
	\underline K_n^{MW}(X,\sL) \hook C^0(X, \underline K^{MW}_n(\sL))=\bigoplus_{x\in X^{(0)}} K_n^{MW}(\kappa(x),\sL)
	\]
	is natural with respect to flat morphisms. This is obvious because the flat pullback on $C^0$ is by definition the sum of the pullbacks in Milnor–Witt $K$-theory.
\end{proof}

If $Z\subset X$ is smooth of codimension $c$, comparing Rost–Schmid complexes yields an isomorphism
\[
\Pi\colon \Chw^n_Z(X,\sL) \simeq \Chw^{n-c}(Z,\sL\otimes \det(\sN_{Z/X})^{-1}),
\]
called the \emph{purity isomorphism}. 

\begin{lem}\label{lem:purity-comparison}
	Under the identification~\eqref{eqn:hzmw-chw}, the purity isomorphism $\Pi$ coincides with the Morel–Voevodsky purity isomorphism~\eqref{eqn:purity-support}.
\end{lem}

\begin{proof}
	We can reduce to the case of the zero section of a vector bundle using the functoriality of the Rost–Schmid complex for smooth morphisms (Lemma~\ref{lem:flat pullback}), Jouanolou devices, and étale neighborhoods (cf.\ \cite[Lemma 3.22]{hoyois-sixops}).
	Thus let $V=\bV(\sE)$ be a vector bundle over $X\in\Sm_k$. We must show that the following square commutes:
	\[
	\begin{tikzcd}
		\pi_0\hzmw_X(V,\xi) \ar{r}{\eqref{eqn:purity-support}} \ar{d}{\simeq} & \pi_0\hzmw(X,\xi-\sE) \ar{d}{\simeq} \\
		\Chw^n_X(V,\det\xi) \ar{r}{\Pi} & \Chw^{n-c}(X,\det(\xi)\otimes \det(\sE)^{-1}).
	\end{tikzcd}
	\]
	The top horizontal map is now the identity map by \cite[Lemma 2.2]{VV}. Each vertical map is the composition of the Thom isomorphism for $\hzmw$ and the canonical map $\hzmw\to \underline K_*^{MW}$.
	Levine shows in \cite[Proposition 3.7]{LevineQuad} that the purity isomorphism $\Pi$ above is the Thom isomorphism of an $\SL$-orientation on the cohomology theory represented by the motivic spectrum $\underline K_*^{MW}$. By \cite[Theorem 5.9]{PaninWalter}, such an orientation is classified by a \emph{unital} morphism of spectra $\mathrm{MSL}\to \underline K_*^{MW}$. But since the unit map $\1\to\mathrm{MSL}$ is a $\underline\pi_0$-isomorphism, there is a unique such morphism. Therefore the map $\hzmw\to \underline K_*^{MW}$ intertwines the respective Thom isomorphisms, which implies the commutativity of the above square.
\end{proof}

\sssec{Comparison of pushforwards in Chow–Witt theory} 
Let $k$ be a perfect field, $f\colon X\to Y$ a morphism between smooth $k$-schemes, $Z\subset X$ a closed subscheme such that the restriction of $f$ to $Z$ is finite, and $\sL$ an invertible sheaf on $Y$. 
We recall the definition of the \emph{Calmès–Fasel pushforward}
\[
f_*\colon \Chw^{n+d}_Z(X,f^*\sL \otimes \omega_f) \to \Chw^{n}_{f(Z)}(Y,\sL),
\]
where $d=\rk(\sL_f)$ \cite[page 10]{Calmes:2014ab}. It is induced by the morphism of Rost–Schmid complexes
\[
f_*\colon C^{*+d}_Z(X,\underline K_{n+d}^{MW}(f^*\sL \otimes \omega_f)) \to C^{*}_{f(Z)}(Y,\underline K_{n}^{MW}(\sL)),
\]
which in degree $i$ is the sum of the absolute transfers \cite[Definition 5.4]{Morel}
\[
\bigoplus_{x\in X^{(i+d)}\cap Z} K_{n-i}^{MW}(\kappa(x), f^*\sL \otimes \omega_f \otimes\omega_x)
\to
\bigoplus_{y\in Y^{(i)}\cap f(Z)} K_{n-i}^{MW}(\kappa(y), \sL \otimes\omega_y).
\]

\begin{prop}\label{prop:pushforward-comparison}
	Let $f\colon X\to Y$ be a morphism between smooth $k$-schemes, let $Z\subset X$ be a closed subscheme such that the restriction of $f$ to $Z$ is finite, and let $\xi\in K(Y)$ be of rank $c=\codim(f(Z),Y)$. Then, under the identification~\eqref{eqn:hzmw-chw}, the Gysin map
	\[
	f_!\colon \hzmw_Z(X,f^*\xi+\sL_f) \to \hzmw_{f(Z)}(Y,\xi)
	\]
	(see \sssecref{sssec:pushforwad-coh}) agrees with the Calmès–Fasel pushforward
	\[
	f_*\colon \Chw{}^{c+d}_Z(X,\det(f^*\xi)\otimes \omega_f) \to \Chw{}^{c}_{f(Z)}(Y,\det(\xi)), \quad d=\rk(\sL_f).
	\]
\end{prop}

\begin{proof}
	Since $c+d$ is the codimension of $Z$ in $X$, we have an exact sequence
	\[
	0 \to \Chw^{c+d}_Z(X,\sL) \to \bigoplus_{x\in X^{(c+d)}\cap Z} GW(\kappa(x),\sL\otimes\omega_x)\xrightarrow{\partial} \bigoplus_{x\in X^{(c+d+1)}\cap Z} W(\kappa(x),\sL\otimes\omega_x),
	\]
	where $\omega_{x}=\omega_{\kappa(x)/\sO_{X,x}}$. For $x\in X^{(c+d)}\cap Z$, the map
	\[
	 \Chw^{c+d}_Z(X,\sL) \to GW(\kappa(x),\sL\otimes\omega_x)
	\]
	is the filtered colimit of the restriction maps
	\[
	\Chw^{c+d}_Z(X,\sL) \to \Chw^{c+d}_{Z\cap U}(U,\sL) \stackrel\Pi\simeq \Chw^0(Z\cap U, \sL\otimes\omega_{Z\cap U/U}),
	\]
	where $U$ ranges over the open subschemes of $X$ containing $x$ and such that $U\cap Z$ is smooth, and $\Pi$ is the purity isomorphism. By Lemma~\ref{lem:purity-comparison} and the fact that Gysin maps are compatible with base change and with the purity isomorphism, we are reduced to the following claim: for $L/K$ a finite extension of finitely generated fields over $k$, the Gysin map $\hzmw(\Spec L,\sL_{L/K})\to \hzmw(\Spec K)$ coincides with the absolute transfer $\GW(L,\omega_{L/K})\to \GW(K)$. Without loss of generality, we can assume that $L=K(a)$ for some $a\in L$. Both transfers can then be computed in terms of the factorization 
	 \[
	 \Spec L\stackrel a\hookrightarrow \P^1_K \xrightarrow{p} \Spec K.
	 \]
	 More precisely, we use the geometric description of the absolute transfer from \cite[page 99]{Morel}. It is the lower composition in the following diagram:
\[
\begin{tikzcd}
	 \hzmw(\Spec L,\sL_{L/K}) \ar{d}{\Pi}[swap]{\simeq} \ar{r}{a_!} & \hzmw(\P^1_K,\sL_{\P^1_K/K}) & \hzmw(\Spec K) \ar{d}{\Pi}[swap]{\simeq} \ar{l}[swap]{s_!} \\
	 \hzmw_{\Spec L}(\P^1_K,\sL_{\P^1_K/K}) \ar{r} & \hzmw_{\P^1_K-\infty}(\P^1_K,\sL_{\P^1_K/K}) \ar{u} & \hzmw_0(\P^1_K,\sL_{\P^1_K/K}), \ar{l}{\simeq}
\end{tikzcd}
\]
where $s\colon \Spec K \hookrightarrow \P^1_K$ is the zero section. The commutativity of each square is an instance of \sssecref{sssec:pushforwad-coh}(1). Since $p_!s_!=\id$, it follows that the absolute transfer coincides with the Gysin map $p_!a_!$.
\end{proof}

\sssec{$\hzmw$-correspondences vs.\ Milnor–Witt correspondences} 
Let $S$ be essentially smooth over a field. By $\cormw_S$ we denote Calmès and Fasel's category of finite Milnor–Witt correspondences between smooth separated $S$-schemes, as defined in \cite{Calmes:2014ab}.

\begin{prop}
	Let $S$ be essentially smooth over a field. Then the symmetric monoidal $\infty$-category $\Span^{\hzmw}(\Sm_S)$ is a $1$-category and is equivalent to $\cormw_S$.
\end{prop}

\begin{proof}
	Since $\hzmw$ is in the heart of the effective homotopy $t$-structure, $\Span^{\hzmw}(\Sm_S)$ is a $1$-category by Proposition~\ref{prop:discrete}.
	For smooth $S$-schemes $X$ and $Y$, we have
		$$
		\Corr^{\hzmw}_S(X, Y) \simeq  \colim_{Z \subset X \times_S Y} \hzmw_Z(X \times_S Y, \pi_Y^* \Omega_{Y/S})
		$$
		by~\sssecref{sssec:borel-moore vs support}.
	It follows from~\eqref{eqn:hzmw-chw} that 
	\[
	\pi_0\hzmw_Z(X \times_S Y, \pi_Y^* \Omega_{Y/S})\simeq \Chw^d_Z(X\times_SY,\pi_Y^*\omega_{Y/S}).
	\] 
	In particular, the mapping spaces in $\Span^{\hzmw}(\Sm_S)$ are the same as in $\cormw_S$.
	
	To compare the composition laws, we use the description of the composition in $\Span^{\hzmw}(\Sm_S)$ via the pairing $\theta$ (Lemma~\ref{lem:compositions-agree}). The composition in $\cormw_S$ is defined in the same way, except that it uses the Calmès–Fasel pushforward instead of the Gysin map, but these are the same by Proposition~\ref{prop:pushforward-comparison}. Finally, the symmetric monoidal structures agree by the multiplicativity of the isomorphisms~\eqref{eqn:hzmw-chw}.
\end{proof}

\sssec{Compatibility with the functor of D\'eglise–Fasel}
\label{sssec: D-F functor}
Let $k$ be a perfect field. In \cite[Proposition 2.1.12]{DegliseFasel}, Déglise and Fasel define a functor
$$\alpha \colon \Span^\efr_*(\Sm_k) \to \cormw_k, $$
where $\Span^\efr_*(\Sm_k)$ is the category whose objects are smooth $k$-schemes (separated and of finite type) and whose mapping spaces are the sets
\[
\bigvee_{n\geq 0}\Corr^{\efr,n}_k(X,Y)
\]
(see \cite[3.4.7]{EHKSY1}). The rest of this section will be devoted to the proof of the following comparison theorem:

\begin{thm}
	\label{thm: comparison of functors}
	Let $\lambda \colon \Span^\efr_*(\Sm_k) \to \ho\Span^\fr(\Sm_k)$ be the functor defined in~\cite[3.4.7]{EHKSY1}. Then the following diagram commutes:
	\begin{equation*}
	\begin{tikzcd}
	\Span^\efr_*(\Sm_k) \ar{d}[swap]{\lambda} \ar{r}{\alpha} & \cormw_k.  \\
	\ho\Span^\fr(\Sm_k) \ar{ru}[swap]{\Phi^{\hzmw}} &
	\end{tikzcd}
	\end{equation*}
\end{thm}

\sssec{}
We briefly recall the construction of the functor $\alpha$. On objects one has $\alpha(X) = X$. 
Given an equationally framed correspondence $c = (Z, U, \phi, g)\in\Corr^{\efr,n}_k(X,Y)$, we construct a finite MW-correspondence $\alpha(c) \in \cormw_k(X, Y)$ as follows. 

Write $\phi = (\phi_1, \dots, \phi_n)$ and denote by $|\phi_i|$ the vanishing locus of $\phi_i\colon U\to \A^1$, so that $Z = |\phi_1| \cap \dots \cap |\phi_n|$. Since $Z\subset U$ is everywhere of codimension $n$, $|\phi_i|$ does not contain any generic point of $U$. 
Each $\phi_i$ can thus be seen as an element of $\bigoplus_{u\in U^{(0)}} \kappa(u)^\times$ and hence defines an element $[\phi_i]\in \bigoplus_{u\in U^{(0)}}\K_1^{MW}(\kappa(u))$. Applying the residue map
$$ \partial \colon \bigoplus_{u\in U^{(0)}}\K_1^{MW}(\kappa(u)) \longrightarrow \bigoplus_{x \in U^{(1)}} \K_0^{MW}(\kappa(x), \omega_x)$$
we obtain an element $\partial[\phi_i]$ supported on $|\phi_i|$, which defines a cohomology class
$$Z(\phi_i) \in H^1_{|\phi_i|}(U, \underline K_1^{MW}).$$ 
Using the product in Milnor–Witt K-theory, we get an element
$$Z(\phi) = Z(\phi_1) \cdot \ldots \cdot Z(\phi_n) \in H^n_Z(U, \underline K_n^{MW})= \Chw_Z^n(U).$$

The \'etale morphism $u \colon U \to \A^n_X$ induces a trivialization $\omega_{U/X} \simeq u^* \omega_{\A^n_X/X} \simeq \sO_U$. 
Denote by $\pi \colon \A^n_X \to X$ the projection.
Since $Z$ is finite and equidimensional over $X$, the morphism $(\pi u, g) \colon U \to X \times Y$ sends $Z$ to a closed subscheme $T$, which is finite and equidimensional over $X$. 
The finite MW-correspondence $\alpha(c) \in \cormw_k(X, Y)$ is then the image of $Z(\phi)$ by the Calmès–Fasel pushforward
\[
(\pi u, g)_* \colon  \Chw^n_Z(U)\simeq \Chw^n_Z(U,\omega_{U/X}) \longrightarrow \Chw^{d}_T(X \times Y,\omega_{X \times Y / X}),
\quad d=\dim(Y).
\]

\sssec{} 

The first step in the proof of Theorem~\ref{thm: comparison of functors} is to recast the construction of $Z(\phi)$ as a Thom class.

We recall that for $E$ a motivic ring spectrum, the Thom class of a locally free sheaf $\sE$ on $X$ is the image of $1$ by the purity equivalence $E(X) \simeq E_X(\bV(\sE),\sE)$. By Lemma~\ref{lem:purity-comparison}, the Thom class of $\sE$ in Chow–Witt theory has an explicit representative in the Rost–Schmid complex, namely
\[
1\in \bigoplus_{x\in X^{(0)}} \GW(\kappa(x)) \simeq C^{n}_{X}(\bV(\sE), \underline K_n^{MW}(\det\sE)),\quad n=\rk\sE.
\]

\begin{lem}
\label{lem: thom}
Let  $X, Y \in \Sm_k$ and let $c = (Z, U, \phi, g) \in \Corr^{\efr,n}_k(X, Y)$. Assume that the morphism $\phi \colon U \to \A^n_k$ is flat. Let $t_n \in \Chw_0^n(\A^n_k)$ be the Thom class of the trivial vector bundle $\A^n_k\to \Spec k$. 
Then $Z(\phi) = \phi^* (t_n)$ in $\Chw^n_Z(U)$.
\end{lem}

\begin{proof}
It is enough to show that $Z(\phi_i) = \phi_i^* (t_1)$, because $Z(\phi) = Z(\phi_1) \cdot \ldots \cdot Z(\phi_n)$ and the Thom class is multiplicative with respect to direct sum of vector bundles. 
By Lemma~\ref{lem:flat pullback}, since $\phi_i\colon U\to\A^1_k$ is flat, the pullback $\phi_i^*$ on Chow–Witt groups can be computed using Rost–Schmid complexes. The commutative square
\[
\begin{tikzcd}
	C^0(\A^1_k,\underline K_1^{MW}) \ar{d}[swap]{\partial} \ar{r}{\phi_i^*} & C^0(U,\underline K_1^{MW}) \ar{d}{\partial} \\
	C^1(\A^1_k,\underline K_1^{MW}) \ar{r}{\phi_i^*} & C^1(U,\underline K_1^{MW})
\end{tikzcd}
\]
shows that $\partial[\phi_i]=\phi_i^*(\partial[\id_{\A^1}])$ in $C^1_{|\phi_i|}(U,\underline K_1^{MW})$, hence that $Z(\phi_i)=\phi_i^*(Z(\id_{\A^1}))$ in $\Chw_{|\phi_i|}^1(U)$.
It remains to observe that $Z(\id_{\A^1})=t_1$, since the residue map $\partial_t\colon K_1^{MW}(k(t)) \to GW(k)$ takes $[t]$ to $1$.
\end{proof}

\sssec{}
The following lemma shows that the flatness assumption in Lemma~\ref{lem: thom} is essentially vacuous.

\begin{lem}
\label{lemma: flat}
Let $S$ be a regular Noetherian scheme and $X,Y\in \Sch_S$.
Suppose that $X$ is flat over $S$.
Then, for every $(Z, U, \phi, g) \in  \Corr^{\efr,n}_S(X, Y)$, the morphism $\phi\colon U\to \A^n_S$ is flat in an open neighborhood of $Z$. 
\end{lem}
\begin{proof}
It suffices to show that $\phi\colon U\to \A^n_S$ is flat at every point $z\in Z$.
Let $s=\phi(z)$ and let $(x_1,\dotsc,x_d)\in \Oo_{S,s}$ be a regular system of parameters.
By \cite[Proposition 2.1.18]{EHKSY1}, $Z$ is flat over $S$ and $(\phi_1,\dotsc,\phi_n)$ is regular sequence in $\Oo_{U,z}$ with quotient $\Oo_{Z,z}$. The flatness of $Z$ implies that the image of $(x_1,\dotsc,x_d)$ in $\Oo_{Z,z}$ is a regular sequence. Thus, the local homomorphism $\phi^* \colon \Oo_{\A^n_S, s} \to \Oo_{U, z}$ sends the regular system of parameters $(t_1,\dotsc,t_n,x_1,\dotsc,x_d)$ to a regular sequence. The criterion of \cite[Tag 07DY]{stacks} now shows that $\phi$ is flat.
\end{proof}

\begin{proof}[Proof of Theorem~\ref{thm: comparison of functors}]
Let $X,Y\in \Sm_k$ and let $c = (Z, U, \phi, g) \in \Corr^{\efr,n}_k(X,Y)$ be an equationally framed correspondence:
\[
\begin{tikzcd}
  & U \ar{dr}{(\phi, g)} \ar{dl}[swap]{u} &  \\
  \A^n_X \ar{d}[swap]{\pi} & Z \ar[hook]{u}{i} \ar{dl}{f} \ar[hook']{l} \ar{dr}[swap]{h} & \A^n \times Y \\ 
 X & & Y. \ar[hook]{u}[swap]{0}
\end{tikzcd}
\]
Let $\tau\in \Maps_{K(Z)}(0,\sL_f)$ be the induced trivialization of the cotangent complex and let $T=(f,h)(Z)\subset X\times Y$. 
By Lemma~\ref{lemma: flat} we can assume that $\phi\colon U\to \A^n_k$ is flat. 
Combining \sssecref{sssec:bm-support}(2) and Proposition~\ref{prop:pushforward-comparison}, we have a commutative square
\[
\begin{tikzcd}
	\hzmw^\BM(Z/X) \ar{r}{(f,h)_*} \ar{d}[swap]{\simeq} & \hzmw^\BM(T/X) \ar{d}{\simeq} \\
	\Chw^n_Z(U) \ar{r}{(\pi u,g)_*} & \Chw^d_T(X\times Y,\omega_{X\times Y/X}).
\end{tikzcd}
\]
It is thus enough to show that $Z(\phi)\in \Chw^n_Z(U)$ corresponds to $\tau_*(\eta_f) \in \hzmw^\BM(Z/X)$ under the left-hand isomorphism. We have
\begin{align*}
	Z(\phi) & =\phi^*(t_n) & \text{(Lemma~\ref{lem: thom})} \\
	& = \phi^*0_!(1) & \text{(definition of the Thom class)} \\
	& = \tau_*i_!(1), & \text{(base change)}
\end{align*}
where $0_!$ and $i_!$ are the Gysin maps defined in~\sssecref{sssec:pushforwad-coh}. To conclude, note that the isomorphism
\[
\hzmw_Z(U,\sN_i) \simeq \hzmw^\BM(Z/X,\sL_f)
\]
sends $i_!(1)$ to $\eta_f$, by definition of $i_!$.
\end{proof}


\let\mathbb=\mathbf

{\small
\newcommand{\etalchar}[1]{$^{#1}$}
\providecommand{\bysame}{\leavevmode\hbox to3em{\hrulefill}\thinspace}

}

\parskip 0pt

\end{document}